\documentclass[Svgc, numbook, envcountsame]{Svgc}
\usepackage{latexsym}
\usepackage{graphics}
\usepackage{amsmath, amssymb}

\def\beqa{\begin{eqnarray*}}
\def\eeqa{\end{eqnarray*}}
\newcommand{\abs}[1]{\left\vert#1\right\vert}

\journalname{Graphs and Combinatorics}

\begin{document}
\title{Monochromatic and Zero-Sum Sets of Nondecreasing Diameter}

\author{Carl R. Yerger \inst{1} }

\thanks{The author was supported by NSF grant DMS0097317.  $Email$: cyerger@math.gatech.edu}


%
%
\institute{Georgia Institute of Technology}
\maketitle
\begin{abstract}
Let $k, r, s \in \mathbb{N}$ where $r \geq s \geq 2$.  Define
$f(s,r,k)$ to be the smallest positive integer $n$ such that for
every coloring of the integers in $[1,n]$ there exist subsets
$S_1$ and $S_2$ such that:  $(a)$ $S_1$ and $S_2$ are
monochromatic (but not necessarily of the same color), $(b)$
$|S_1| = s$, $|S_2| = r, $ $(c)$ $\max(S_1) < \min(S_2)$, and
$(d)$ $diam(S_1) \leq diam(S_2)$.  In this paper, we prove the
following:
\begin{center} For $r \geq s \geq 2$, $f(s,r,2) = $
\end{center}
\begin{displaymath}
\left \{\begin{array}{ll} 5s - 3
& \textrm{ if $s = r$ } \\
4s + r - 3 & \textrm{ if $ s < r \leq 2s - 2$ } \\
2s + 2r - 2 & \textrm{ if $r > 2s - 2$  } \end{array} \right.
\end{displaymath}
\begin{center} For $r \geq s \geq 3$, $f(s,r,\mathbb{Z}) = $
\end{center}
\begin{displaymath}
\left \{\begin{array}{ll} 5s - 3
& \textrm{ if $s = r$ } \\
4s + \max(r, s + \frac{s}{(r,s)} - 1) - 3 & \textrm{ if $ s < r \leq 2s - 2$ } \\
2s + 2r - 2 & \textrm{ if $r > 2s - 2$  } \end{array} \right.
\end{displaymath}
\begin{center} For $r \geq s \geq 3$, $f(s,r,3) = $
\end{center}
\begin{displaymath}
\left \{\begin{array}{ll} 9s - 7
& \textrm{ if $r \leq \lfloor \frac{5s - 1}{3} \rfloor - 1$ } \\
4s + 3r - 4 & \textrm{ if $ \lfloor \frac{5s - 1}{3} \rfloor - 1 < r \leq 2s - 2$ } \\
6s + 2r - 6 & \textrm{ if $ 2s - 2  < r \leq 3s - 3$ } \\
3s + 3r - 4 & \textrm{ if $r > 3s - 3$ ($r \geq s \geq 2$) }
\end{array} \right.
\end{displaymath}
\begin{center} For $r \geq s \geq 3$, $f(s,r,\{\infty\} \cup \mathbb{Z}) = $
\end{center}
\begin{displaymath}
\left \{\begin{array}{ll} 9s - 7 + \frac{s}{(r,s)} - 1
& \textrm{ if $r \leq \lfloor \frac{5s - 1}{3} \rfloor - 1$ } \\
3s + 2r - 3 + \min( 3s - 3, s + r - 1 + \frac{s}{(r,s)} - 1) & \textrm{ if $ \lfloor \frac{5s - 1}{3} \rfloor - 1 < r \leq 3s - 3$ } \\
3s + 3r - 4 & \textrm{ if $r > 3s - 3$ ($r \geq s \geq 2$) }
\end{array} \right.
\end{displaymath}

We prove that the theorems defining $f(s,r,2)$ and $f(s,r,3)$
admit a partial generalization in the sense of the
Erd\H{o}s-Ginzburg-Ziv theorem.  This work begins the off-diagonal
case of the results of Bialostocki, Erd\H{o}s, and Lefmann.

\end{abstract}
\receive{December 2005} \finalreceive{}
\section{Introduction}
In 1961, Erd\H{o}s, Ginzburg and Ziv in \cite{egz} proved the
following theorem that has been the subject of many recent
developments.

\begin{theorem}
Let $m \in \mathbb{N}$.  Every sequence of $2m - 1$ elements from
$\mathbb{Z}$ contains a subsequence of $m$ elements whose sum is
zero modulo $m$.
\end{theorem}

Notice that this theorem is a generalization of the pigeonhole
principle.  Indeed, if the sequence contains only the residues 0
and 1, then this theorem describes the situation of placing $2m -
1$ pigeons in 2 holes.

We begin by introducing some notation and definitions.   A mapping
$\Delta : X \longrightarrow C$, where $C$ is the set of colors, is
called a $coloring$.  If $C = \{1,\ldots, k\}$, we say that
$\Delta$ is a $k$-coloring.  If the set of colors is $\mathbb{Z}$,
the additive group of the integers, we call $\Delta$ a
$\mathbb{Z}$-coloring. A set $X$ is called $monochromatic$ if
$\Delta(x) = \Delta(x')$ for all $x, x' \in X$. In a
$\mathbb{Z}$-$coloring$ of $X$, a subset $Y$ of $X$ is called
$zero$-$sum$ mod $m$ if $\sum_{y \in Y} \Delta(y)$ $\equiv 0$ mod
$m$. If $a, b \in \mathbb{N}$, then $[a,b]$ is the set of integers
$\{n \in \mathbb{N} | a \leq n \leq b\}$.  For finite $X \subseteq
\mathbb{N}$, define the $diameter$ of $X$, denoted by $diam(X)$,
to be $diam(X) = \max(X) - \min(X)$.  For integers, $s$, $r$, let
$(s,r)$ be the greatest common factor of $s$ and $r$.  Denote by
$\infty$ an additional color, which does not belong to
$\mathbb{Z}$. For brevity and ease of expression, we denote proofs
of a parallel nature by parentheses.
\newline
\begin{definition}
Let $f(s,r,k)$ $($Let $f(s, r, \mathbb{Z}))$ $($Let $f(s, r,
\{\infty\})$ $)$ be the the smallest positive integer $n$ such
that for every coloring $\Delta : [1, n] \longrightarrow [1, k]$ $
(\Delta : [1, n] \longrightarrow \mathbb{Z})$ $ (\Delta : [1, n]
\longrightarrow \{\infty\} \cup \mathbb{Z})$, there exist two
subsets $S_1, S_2$ of $[1, n]$, which satisfy:
\newline
$(a) S_1 $ and $S_2$ are monochromatic $(S_1 $ is zero-sum mod $s$
and $S_2$ is zero-sum mod $r)$  $(S_1 $ is either
$\infty$-monochromatic or zero-sum mod $s$ and $S_2$ is either
$\infty$-monochromatic or zero-sum mod $r)$,
\newline $(b)$$ |S_1| = s, |S_2| = r$, \newline $(c) \max(S_1) <
\min(S_2)$, and \newline $(d)$$ diam(S_1) \leq diam(S_2)$.
\end{definition}

In \cite{bial} it was shown that $f(m,m,2) = f(m,m,\mathbb{Z}) =
5m - 3$ and $f(m,m,3) = f(m,m,$ $\{\infty\} \cup \mathbb{Z} ) = 9m
- 7$. (Such theorems are known as zero-sum generalizations in the
sense of Erd\H{o}s-Ginzberg-Ziv).  Several papers continued the
investigations of \cite{bial}.  The body of work done on this
topic is mostly separated into two sections. The first involves
coloring the integers with one set, begun in \cite{schaal} and
continued in \cite{guy}, \cite{rash5}, \cite{rash2} and
\cite{rash3}.  The second involves coloring the integers with two
sets, begun in \cite{bial} and further considered in \cite{rash1},
\cite{davidandy}, \cite{schultz} and \cite{biji}.  Initial work on
coloring the integers with three sets has begun in
\cite{davidcarl}.

In this paper, we evaluate $f(s,r,2),$ $f(s,r,\mathbb{Z})$,
$f(s,r,3)$ and $f(s,r,\{\infty\} \cup \mathbb{Z})$ for one of the
off-diagonal cases, namely $r \geq s$, as shown in the tables
below.

\begin{center} For $r \geq s \geq 2$, $f(s,r,2) = $
\end{center}
\begin{displaymath}
\left \{\begin{array}{ll} 5s - 3
& \textrm{ if $s = r$ } \\
4s + r - 3 & \textrm{ if $ s < r \leq 2s - 2$ } \\
2s + 2r - 2 & \textrm{ if $r > 2s - 2$  } \end{array} \right.
\end{displaymath}
\begin{center} For $r \geq s \geq 3$, $f(s,r,\mathbb{Z}) = $
\end{center}
\begin{displaymath}
\left \{\begin{array}{ll} 5s - 3
& \textrm{ if $s = r$ } \\
4s + \max(r, s + \frac{s}{(r,s)} - 1) - 3 & \textrm{ if $ s < r \leq 2s - 2$ } \\
2s + 2r - 2 & \textrm{ if $r > 2s - 2$  } \end{array} \right.
\end{displaymath}
\begin{center} For $r \geq s \geq 3$, $f(s,r,3) = $
\end{center}
\begin{displaymath}
\left \{\begin{array}{ll} 9s - 7
& \textrm{ if $r \leq \lfloor \frac{5s - 1}{3} \rfloor - 1$ } \\
4s + 3r - 4 & \textrm{ if $ \lfloor \frac{5s - 1}{3} \rfloor - 1 < r \leq 2s - 2$ } \\
6s + 2r - 6 & \textrm{ if $ 2s - 2  < r \leq 3s - 3$ } \\
3s + 3r - 4 & \textrm{ if $r > 3s - 3$ ($r \geq s \geq 2$) }
\end{array} \right.
\end{displaymath}
\begin{center} For $r \geq s \geq 3$, $f(s,r,\{\infty\} \cup \mathbb{Z}) = $
\end{center}
\begin{displaymath}
\left \{\begin{array}{ll} 9s - 7 + \frac{s}{(r,s)} - 1
& \textrm{ if $r \leq \lfloor \frac{5s - 1}{3} \rfloor - 1$ } \\
3s + 2r - 3 + \min( 3s - 3, s + r - 1 + \frac{s}{(r,s)} - 1) & \textrm{ if $ \lfloor \frac{5s - 1}{3} \rfloor - 1 < r \leq 3s - 3$ } \\
3s + 3r - 4 & \textrm{ if $r > 3s - 3$ ($r \geq s \geq 2$) }
\end{array} \right.
\end{displaymath}

A collection of sets $S_i$ that satisfy conditions $(a)$, $(b)$,
$(c)$ and $(d)$ of Definition 1.2 is called a $solution$ to
$f(s,r,k)$ ($ f(s,r,\mathbb{Z})$),  ($f(s,r,\{\infty\} \cup
\mathbb{Z})$).  For every positive integer $r$ and $s$, it turns
out that $f(s,r,2)$ and $f(s,r,3)$ admit a zero-sum generalization
in the sense of Erd\H{o}s-Ginzburg-Ziv only for some values of $s$
and $r$.  This paper is organized as follows.  In Section 2, we
present some preliminaries from additive number theory.  In
Section 3, we prove the values of $f(s,r,2)$ and
$f(s,r,\mathbb{Z})$.  In Section 4, the values of $f(s,r,3)$ and
$f(s,r,\{\infty\} \cup \mathbb{Z})$ are determined.

\section{Preliminaries from Additive Number Theory}

The following generalization of the Erd\H{o}s-Ginzburg-Ziv (EGZ)
theorem was shown in \cite{residues3}, and it will be used in our
evaluations of $f(s, r, \mathbb{Z})$, $f(s, r, 3)$ and $f(s, r,
 \{\infty\} \cup \mathbb{Z})$.
\begin{theorem}
Suppose that $S$ is a set.  Let $m \geq 3$ be an integer, let
$\mathbb{Z}_m$ denote the cyclic group on $m$ elements, and let
$\Delta : S \rightarrow \mathbb{Z}_m$ be a coloring for which
$|\Delta(S)| \geq 3$.  If $|S|  = 2m - 2,$ then there exist
distinct integers $x_1, \ldots, x_m$ such that $\sum_{i-1}^m
\Delta(x_i) = 0.$
\end{theorem}
Next, we state a slightly stronger form of the EGZ theorem
\cite{carosurvey} used in our determination of $f(s, r, \{\infty\}
\cup \mathbb{Z}).$
\begin{theorem}
Let $k$ and $m$ be positive integers such that $k|m$.  If
$\Delta:[1, m + k - 1] \rightarrow \mathbb{Z}_m$, then there exist
distinct integers  $x_1, \ldots, x_m \in [1, m+k-1]$ such that
$\sum_{i=1}^m \Delta(x_i) \equiv 0 \, \text{mod} \, k$. Moreover,
$m + k - 1$ is the smallest number for which the above assertion
holds.
\end{theorem}
Throughout the paper, we apply Theorem 2.2 to the following
situation: Suppose there are $a$ elements all colored by the same
residue mod $r$.  This means there are $\frac{s}{(r,s)}$ possible
residues mod $s$ for these $a$ elements.  Therefore, by Theorem
2.2, it takes at most $s + \frac{s}{(r,s)} - 1$ elements to obtain
an $s$-element zero-sum mod $s$ set from the $a$ elements.

Let $G$ be an finite abelian group.  If $a_1, \ldots, a_s$ is a
sequence $S$ of elements from $G$, then an $n$-$set$ $partition$
of $S$ is a collection of $n$ nonempty subsequences of $S$,
pairwise disjoint as sequences, such that every term of $S$
belongs to exactly one of the subsequences, and the terms in each
subsequence are all distinct, so that the subsequences may be
considered sets.  Grynkiewicz proved in \cite{davidnumbertheory}
and \cite{rash1} the following theorem, which is used in the
determination of $f(s, r, \{\infty\} \cup \mathbb{Z})$.   The
result is an extension of the Cauchy-Davenport theorem.
\begin{theorem}
Let $a_1,\ldots,a_s$ be a sequence $S$ of elements from an abelian
group $G$ of order $m$ with an $n$-set partition
$P=P_1,\ldots,P_n$. Furthermore, let $p$ be the smallest prime
divisor of $m$. Then either: \noindent \newline(i) there exists an
$n$-set partition $A=A_1,A_2,\ldots,A_n$ of $S$ such that: \beqa
\abs{\sum_{i=1}^n A_i}\geq
\min\left\{m,\,(n+1)p,\,|S|-n+1\right\};\eeqa  furthermore, if
$n'\geq\frac{m}{p}-1$ is an integer such that $P$ has at least
$n-n'$ cardinality one sets and if $|S|\geq n+\frac{m}{p}+p-3$,
then we may assume there are at least $n-n'$ cardinality one sets
in $A$, or $(ii)$: \newline $(a)$ there exists $\alpha\in G$ and a
nontrivial proper subgroup $H_a$ of index $a$ such that all but at
most $a-2$ terms of $S$ are from the coset $\alpha +H_a$; and
\noindent \newline  $(b)$ there exists an $n$-set partition
$A_1,A_2,\ldots,A_n$ of the subsequence of $S$ consisting of terms
from $\alpha+H_a$ such that $\sum_{i=1}^n A_i=n\alpha+H_a.$
\end{theorem}

\section{The 2-Coloring and its EGZ Generalization}

In this section, we will evaluate $f(s,r,2)$ and
$f(s,r,\mathbb{Z})$.  First, we introduce notion helpful in the
proof. For $c \in C$, where $\Delta^{-1}(c) \neq \emptyset$, we
denote $first(c) = \min\{x \in X | \Delta(x) = c\}$ and $last(c) =
\max\{x \in X|\Delta(x) = c\}$.  A coloring $\Delta : [1, n]
\longrightarrow C$ will be identified by the strings
$\Delta(1)\Delta(2) \cdots \Delta(n)$, and we use $x^i$ to denote
the string, $xx\cdots x$, of length $i$.  First, we evaluate
$f(s,r,2)$.

\begin{lemma}
Let $r$ and $s$ be positive integers with $r \geq s \geq 2$ and $r
\leq 2s - 2$.  If $\Delta : [1, 2s + r - 2]$ $\longrightarrow$ \
$\mathbb{Z}$ is a coloring, then either,
\newline
(i) there exists an $r$-element zero-sum mod $r$ subset $R$
$\subseteq$ $[1, 2s + r - 2]$ such that $diam(R) \geq 2s - 2$ or
\newline (ii) there exists a zero-sum solution to $f(s, r,
\mathbb{Z})$.
\end{lemma}
\begin{proof}
Let $I_1 = [1,r], I_2 = [r + 1, 2s - 2]$, and $I_3 = [2s - 1, 2s +
r - 2]$.  Since $|I_1 \cup (I_3 \setminus \{2s - 1\})| = |(I_1
\setminus \{r\}) \cup I_3| = 2r - 1$, by the pigeonhole principle,
it follows that both $I_1 \cup (I_3 \setminus \{2s - 1\})$ and
$(I_1 \setminus \{r\}) \cup I_3$ contain $r$-element monochromatic
sets, say $R_1$ and $R_2$, respectively.  If $R_1 \cup I_1 \neq
\emptyset$ and $R_1 \cap (I_3 \setminus \{2s - 1\}) \neq
\emptyset$, then $diam(R_1) \geq 2s - 2$ whence $(i)$ follows.
Similarly, if $R_2 \cap (I_1 \setminus \{r\}) \neq \emptyset$ and
$R_2 \cap I_3 \neq \emptyset$, then $diam(R_2) \geq 2s - 2$, and
$(i)$ follows.  Therefore, we can assume that $I_1$ and $I_3$ are
$r$-element monochromatic sets.  Therefore, condition $(ii)$
follows by taking $[1, s]$ and $I_3$. \qed
\end{proof}

\begin{theorem}
Let $r$ and $s$ be positive integers where $r \geq s \geq 2$.  If
$r \leq 2s - 2$, then $f(s,r,2) = 4s + r - 3$.
\end{theorem}
\begin{proof}
The coloring $\Delta : [1, 4s + r - 4] \longrightarrow \{0,1\}$
given by the string:
\begin{align*}
01^{s-1}0^{s-1}1^{2s-2}0^{r-1}
\end{align*}
implies that $f(s,r,2) > 4s + r - 4.$  Next we show that $f(s,r,2)
\leq 4s + r - 3$.  Let $\Delta : [1, 4s + r - 3] \longrightarrow
\{0,1\}$ be an arbitrary coloring.  By the pigeonhole principle,
the interval $[1, 2s - 1]$ contains an $s$-element monochromatic
subset $K$ such that $diam(K) \leq 2s - 2$.  Shifting the interval
$[2s, 4s + r - 3]$ to the interval $[1, 2s + r - 2]$ and applying
Lemma 3.1 completes the proof. \qed
\end{proof}

We now will compute $f(s, r, \mathbb{Z})$.  First, consider a
special case of the $2$-coloring that allows us to give a simple
example of how $f(s,r,2) \neq f(s,r,\mathbb{Z})$.

\begin{theorem}
Let $r$ and $s$ be positive integers where $r \geq s \geq 3$ and
$r,s$ are coprime.  If $r \geq 2s - 2$, then $f(s,r,\mathbb{Z}) =
6s - 4.$
\end{theorem}
\begin{proof}
The coloring $\Delta : [1, 6s - 5] \longrightarrow \{1,2\}$ given
by the string:
\begin{align*}
01^{s-1}0^{s-1}1^{s-1}0^{s-1}1^{s-1}0^{s-1} & \textnormal{   mod } s \\
01^{s-1}0^{s-1}1^{2s - 2}0^{2s - 2} & \textnormal{   mod } r
\end{align*}
implies that $f(s,r,2) > 2s + 2r - 3$.  By Theorem 2.1, it follows
that the interval $[1, 2s - 1]$ contains an $s$-element zero-sum
mod $s$ subset $S_1$ with $diam(S_1) \leq 2s - 2$, and that the
interval $[2s, 6s - 4]$ contains an $r$-element zero-sum mod $r$
subset $S_2$ with $diam(S_2) > 2s - 2.$  Hence, $S_1$ and $S_2$
complete the proof. \qed
\end{proof}

As we can see, since $f(s,r,2) = 4s + r - 3 \neq 6s - 5$ (for $r <
2s - 2$), in this case, we do not have a generalization in the
sense of Erd\H{o}s-Ginzburg-Ziv.  However, when $r > 2s - 2$, an
EGZ generalization exists.

\begin{theorem}  Let $r$ and $s$ be positive integers where $r
\geq s \geq 2$.  If $r > 2s - 2$, then $f(s,r,2) =
f(s,r,\mathbb{Z}) = 2s + 2r - 2$.
\end{theorem}

\begin{proof}  The coloring $\Delta : [1, 2s  + 2r - 3]
\longrightarrow \{1,2\}$ given by the string
\begin{align*}
01^{s-1}0^{s-1}1^{r-1}0^{r-1} & \textnormal{   mod } s \\
01^{s-1}0^{s-1}1^{r-1}0^{r-1} & \textnormal{   mod } r
\end{align*}
implies that $f(s,r,2) > 2s + 2r - 3$.  By Theorem 1.1, it follows
that the interval $[1, 2s - 1]$ contains an $s$-element zero-sum
mod $s$ subset $S_1$ with $diam(S_1) \leq 2s - 2$, and that the
interval $[2s, 2s + 2r - 2]$ contains an $r$-element zero-sum mod
$r$ subset $S_2$ with $diam(S_2) > 2s - 2$.  Hence, $S_1$ and
$S_2$ complete the proof, since $f(s,r,2) \geq f(s, r,
\mathbb{Z})$ holds trivially. \qed
\end{proof}

We expand on the previous definitions of $first(k)$ and $last(k)$
in the following manner. Let $x_1<x_2<\ldots<x_n$ be the integers
colored by $c$ in $Y$. Then, for integers $i$ and $j$ such that
$1\leq i\leq j\leq n$, we use the notation $first_i^j(c,Y)$ to
denote $\{x_i,\,x_{i+1},\,\ldots,\,x_j\}$. Likewise
$first_i(c,Y)=\{x_i\}$, and $first(c, Y)=first_1(c,Y)$. Similarly,
we define
$last_i^j(c,Y)=\{x_{n-i+1},\,x_{n-i},\,\ldots,\,x_{n-j+1}\}$,
$last_i(c,Y)=\{x_{n-i+1}\}$, and $last(c,Y)=last_1(c,Y)$.  If we
do not specify $Y$, then it is assumed that $Y$ is the entire
interval.  We now will compute $f(s,r,\mathbb{Z})$ for the
remaining cases.

\begin{lemma}
Let $r$ and $s$ be positive integers with $r \geq s \geq 3, r \leq
2s - 2$, and $(r,s) > 1$ $($and $s + \frac{s}{(r,s)} - 1 > r)$. If
$\Delta : [1, 2s + r - 2]$ $([1, 3s + \frac{s}{(r,s)} - 3])
\longrightarrow \mathbb{Z}$ is a coloring, then either:
\newline $(i)$ there exists an $r$-element zero-sum mod $r$ subset
$R \subseteq [1, 2s + r - 2]$ $([1, 3s + \frac{s}{(r,s)} - 3])$
such that $diam(R) \geq 2s - 2$ or \newline $(ii)$ there exists a
zero-sum solution to $f(s,r,\mathbb{Z})$.
\end{lemma}
\begin{proof}
Let $I_1 = [1,r], I_2 = [r+1, 2s - 2],$ and $I_3 = [2s - 1, 2s + r
- 2].$ ($I_3 = 2s - 1 \cup last_1^{r-1}([1, 3s + \frac{s}{(r,s)}
-3 ]))$.  Since $|I_1 \cup (I_3 \setminus \{ 2s - 1\})| = |(I_1
\setminus \{ r \} ) \cup I_3 | = 2r - 1$, in view of Theorem 1.1,
it follows that both $I_1 \cup (I_3 \setminus \{2s - 1\})$ and
$(I_1 \setminus \{r\}) \cup I_3$ contain an $r$-element zero-sum
mod $r$ set, say $R_1$ and $R_2$, respectively.  If $R_1 \cap I_1
\neq \emptyset$ and $R_1 \cap (I_3 \setminus \{2s - 1\}) \neq
\emptyset$, then $diam(R_1) \geq 2s - 2$, whence $(i)$ follows.
Similarly, if $R_2 \cap (I_1 \setminus \{r\}) \neq \emptyset$,
then $diam(R_2) \geq 2s - 2$, and $(i)$ follows. Therefore, we can
assume that $I_1$ and $I_3$ are $r$-element zero-sum mod $r$ sets.

If $\Delta(I_1) \cap (\Delta(I_3 \setminus \{2s - 1\})) \neq
\emptyset$, then by replacing an element $x \in I_1$ with some
element $y \in I_3 \setminus \{2s - 1\}$, where $\Delta(x) =
\Delta(y)$, then we obtain a set $R$ with $diam(R) \geq 2s - 2$,
whence $(i)$ follows.  Hence, $\Delta(I_1) \cap \Delta(I_3) =
\emptyset$.  Furthermore, if $I_1$ is monochromatic mod $r$, then
it takes at most $s + \lfloor \frac{s}{(r,s)} \rfloor - 1$
elements to obtain an $s$-element zero-sum mod $s$ set $S_2$.  So
if $r \geq s +  \lfloor \frac{s}{(r,s)} \rfloor  - 1$, $S_1$ and
$I_3$ satisfy condition $(ii)$.

(Otherwise, let $J_2 = \{2s - 1\}$ and let $J_3 =
last_1^{r-1}(I_3)$.  Notice that since $r < s + \lfloor
\frac{s}{(r,s)} \rfloor  - 1$ that $S_2$ does not exist yet in
this case (hence the reason for the extended argument).  First
suppose that there are only two residues colored by $\mathbb{Z}$
mod $r$. There cannot be only one residue mod $r$ colored by
$\mathbb{Z}$ since $\Delta(I_1) \cap \Delta(I_3) = \emptyset$. Let
$A$ be the residue mod $r$ used to color $I_1$, and let $B$ be the
other residue mod $r$.  Notice that $last(A) \leq 2s - 2$ else
condition $(i)$ holds.  Thus, the interval $[2s - 1, 3s +
\frac{s}{(r,s)} - 3]$ is colored entirely by residue $B$ mod $r$,
whence residue $B$ mod $r$ contains at least $r$ elements. Thus,
the interval $[1, s + \frac{s}{(r,s)} - 1]$ is colored entirely by
residue $B$.  Thus, there exists an $s$-element zero-sum mod $s$
set $S_3$ with $diam(S_3) \leq s + \frac{s}{(r,s)} - 1$. Condition
$(ii)$ then holds by taking $S_3$ and $J_2 \cup J_3$. Thus, there
exists at least three residues of $\mathbb{Z}$ mod $r$.

Now we will apply Theorem 2.3 in order to finish this section of
the proof.  Let $A$ and $B$ be the residues mod $r$ as defined in
the previous paragraph.  Let $a = first(\mathbb{Z})$ and $b =
first_2(\mathbb{Z})$.  Let $y = last_2(\mathbb{Z})$ and $z =
last(\mathbb{Z})$. We must satisfy the conditions of Theorem 2.3,
by 1) creating an $(r-2)$-set partition of $[1, 3s +
\frac{s}{(r,s)} - 3]$ with $r-3$ sets of cardinality 2 and one set
of cardinality 3, and 2) fixing two elements that are not part of
the set partition, in this case $a$ and $z$, unless $a$ and $z$
contain the only residues not in $A \cup B$.  In this case, fix
$b$ instead of $a$ or $y$ instead of $z$.  In any case, the
diameter of this set described by Theorem 2.3 is at least $2s -
2$.  Hence, we can apply Theorem 2.3 to this case.

Therefore, there exists an $(r-2)$-set partition of $[1,3s +
\frac{s}{(r,s)} - 3]$ with $m-3$ sets of cardinality 2 and one set
of cardinality 3 as well as two fixed elements. (Notice we have at
least $r+1$ elements not colored by $A$ and thus enough for our
set partition because if the only elements where a third color
occurs was in $I_3$, then we can apply Theorem 2.1, and this case
will be discussed after the argument based on Theorem 2.3.) We
apply Theorem 2.3 to observe that $1)$ either there is an
$r$-element zero-sum mod $r$ set $R_4$ that contains the two fixed
elements, and $r-2$ other elements of the set partition, one from
each set, or $2)$ all but at most $r-2$ of the elements colored by
$\mathbb{Z}$ are colored by elements from the same coset
$(a\mathbb{Z}_m + \alpha)$ of $\mathbb{Z}_m$.  If the first
conclusion occurs, $diam(R_4) \geq 2s - 2,$ whence condition $(i)$
follows.  In the latter case, Theorem 2.3 implies that any subset
of cardinality $(r + \frac{r}{a} - 1 + a - 2) \leq
\lceil\frac{3r}{2}\rceil - 1$ must contain an $r$ element zero-sum
mod $r$ set.  Hence, consider the set $R_5 =
first_1^{\lceil\frac{3r}{4}\rceil}(\mathbb{Z},[1, 3s  +
\frac{s}{(r,s)} - 3]) \cup
last_1^{\lceil\frac{3r}{4}\rceil}(\mathbb{Z},[1, 3s  +
\frac{s}{(r,s)} - 3])$.  Within $R_5$, there exists an $r$-element
zero-sum mod $r$ set $R_6$ with $diam(R_6) \geq 3s - 3$. Condition
$(i)$ holds.)

Thus, $\Delta(I_1)$ is not monochromatic mod $r$.  Let $M = [1,
r-1]$, and let $M' = [2, r]$.  Let $N = [2s, 2s + r - 2] (N =
J_3).$  Since $I_1$ is not monochromatic mod $r$, then either $M$
or $M'$, say $M$, is colored by at least two residues.  Let $P = M
\cup N$.  Since $|\Delta(M)| \geq 2$ and since $\Delta(M) \cap
\Delta(N) = \emptyset $, it follows that $|\Delta(P)| \geq 3$.
Hence, since $|M| = r-1$ and $|N| = r-1$, then by Theorem 2.1, it
follows that there exists an $r$-element zero-sum mod $r$ set $P'
\subset P$ with $P' \cap M \neq \emptyset $ and $P' \cap N \neq
\emptyset$, whence condition $(i)$ follows.

\qed
\end{proof}

\begin{theorem}
Let $r$ and $s$ be positive integers where $r \geq s \geq 3, r
\leq 2s -2$ and $(r,s) > 1$.  If $r \geq s + \frac{s}{(r,s)} - 1$,
then $f(s,r,2) = f(s,r,\mathbb{Z}) = 4s + r - 3.$
\end{theorem}
\begin{proof}

 The coloring $\Delta : [1, 4s  + r - 4]
\longrightarrow \{1,2\}$ given by the string:
\begin{align*}
01^{s-1}0^{s-1}1^{2s-2}0^{r-1} & \textnormal{   mod } s \\
01^{s-1}0^{s-1}1^{2s-2}0^{r-1} & \textnormal{   mod } r
\end{align*}
implies that $f(s,r,\mathbb{Z}) > 4s + r - 4$.  Next we show that
$f(s,r,\mathbb{Z}) \leq 4s + r - 3$.  Let $\Delta : [1, 4s + r -
3] \longrightarrow \mathbb{Z}$ be an arbitrary coloring.  By
Theorem 1.1, it follows that the interval $[1, 2s - 1]$ contains
an $s$-element zero-sum mod $s$ subset $K$ with $diam(K) \leq 2s -
2$.  Shifting the interval $[2s,4s + r - 3]$  to the interval $[1,
2s + r - 2]$ and applying Lemma 3.5 completes the proof, since
$f(s,r,2) \geq f(s, r, \mathbb{Z})$ holds trivially. \qed
\end{proof}

\begin{theorem}
Let $r$ and $s$ be positive integers where $r \geq s \geq 3, r
\leq 2s - 2$ and $(r,s) > 1$.  If $r < s + \frac{s}{(r,s)} - 1$,
then $f(s,r,\mathbb{Z}) = 5s + \frac{s}{(r,s)} - 4.$
\end{theorem}
\begin{proof}
The coloring $\Delta : [1, 5s + \frac{s}{(r,s)} - 5]
\longrightarrow \mathbb{Z}$ given by the string:
\begin{align*}
10^{s-1}1^{s-1}0^{s-1}2^{\frac{s}{(r,s)}-1}0^{s - \frac{s}{(r,s)}}1^{s + \frac{s}{(r,s)} - 2} & \textnormal{   mod } s \\
10^{s-1}1^{s-1}0^{2s-2}1^{s + \frac{s}{(r,s)} - 2} & \textnormal{
mod } r
\end{align*}
implies that $f(s, r, \mathbb{Z}) > 5s + \frac{s}{(r,s)} - 5.$
Next we show that $f(s,r,\mathbb{Z}) \leq 5s + \frac{s}{(r,s)} -
4$.  Let $\Delta : [1, 5s + \frac{s}{(r,s)} - 4] \longrightarrow
\mathbb{Z} $ be an arbitrary coloring.  By Theorem 1.1, the
interval $[1, 2s - 1]$ contains an $s$-element zero-sum mod $s$
subset $K$ such that $diam(K) \leq 2s - 2$.  Shifting the interval
$[2s, 5s + \frac{s}{(r,s)} - 4]$ to the interval $[1, 3s +
\frac{s}{(r,s)} - 3]$ and applying Lemma 3.5 completes the proof.
\qed
\end{proof}

\section{The 3-Coloring and its EGZ Generalization}

Let $\delta = r - (2s - 2)$ for $r \geq 2s-2$.  Otherwise, $\delta
= 0$.  First, we will determine $f(s,r,3)$.

\begin{lemma}
Suppose $r$ and $s$ are positive integers with $r \geq s \geq 3, r
\leq \lfloor \frac{5s-1}{3} \rfloor - 1$ $(\lfloor \frac{5s-1}{3}
\rfloor - 1 \leq r \leq 3s - 3)$ and let $A$ be a subset of $[1,
6s - 5]$ with $|A| \geq 5s - 4$.  $($Let $A$ be a subset of $[1, s
+ 3r - 2 - \delta]$ with $|A| \geq 3r - 1 - \delta)$.  If $\Delta
: \longrightarrow \{1,2\}$ is a coloring then either:
\newline
$(i)$ there exists a monochromatic $r$-element subset $S \subseteq
A$ with $diam(S) \geq 3s - 3$ or
\newline
$(ii)$ there exists a solution to $f(s,r,3)$.
\end{lemma}
\begin{proof}
Let $I_1 = first_1^r(A, [1, 6s - 5])$ $(first_1^r(A, [1, s + 3r -
2 - \delta]))$ and $I_3 = last_1^{\max(2s - 1, r -1)} ($ $ A, [1,
6s - 5])$ $ (I_3 = last_1^{\max(2s - 1, r -1)}(A, [1, s + 3r - 2 -
\delta]))$.  Since $|I_1 \cup I_3| \geq 2r + 1$, it follows that
there exists a monochromatic set $S$ in $I_1 \cup I_3$ with
cardinality $r + 1$.  If $S \cap I_1 \neq \emptyset$ and $S \cap
I_3 \neq \emptyset$, then $diam(S) \geq 3s - 3$, and assertion
$(i)$ follows, since we can construct a monochromatic set $S'$
from $S$ with $|S'| = r$ and $diam(S') \geq 3s - 3$.  Otherwise,
since $|I_1| < r + 1$, it follows that $I_3$ contains $S$.  If
$\Delta(I_1) \cap \Delta(I_3) \neq \emptyset$, we can replace an
element $x \in S$ with a corresponding element $y \in I_1$, where
$\Delta(x) = \Delta(y)$, to construct a set $T$ with $diam(T) \geq
3s - 3$.  Assertion $(i)$ is satisfied, since we can construct a
set $T'$ from $T$ with $|T'| = r$ and $diam(T') \geq 3s - 3$. Thus
assume that $\Delta(I_1) \cap \Delta(I_3) = \emptyset$.  Hence,
$I_1$ and $I_3$ are each monochromatic.  Call an element $x \in
[\min(I_1), \max(I_1)]$ a \textit{hole} of $I_1$ if $x \notin
I_1$.  There are at most $s-1$ holes.  Hence, $diam(I_1) \leq 2s -
2.$  Let $S''$ be any $r$-element subset of $I_3$ with $diam(S'')
\geq diam(I_3)$.  The sets $I_1$ and $S''$ satisfy assertion
$(ii)$, thus completing the proof. \qed
\end{proof}

\begin{lemma}
Let $r, s$ be positive integers where $r \geq s \geq 3$, with $r
\geq \lfloor \frac{5s - 1}{3} \rfloor - 1$ $(\lfloor \frac{5s -
1}{3} \rfloor - 1 < r \leq 3s - 3)$.  If $\Delta : [1, 6s - 5]
([1, s + 3r - 2 - \delta]) \longrightarrow \{1,2,3\}$ is a
coloring, then either:
\newline
$(1)$ there exists a monochromatic $r$-element subset $S \subseteq
[1, 6s - 5]$ $(S \subseteq [1, s + 3r - 2 - \delta])$ with
$diam(S) \geq 3s - 3$ or \newline $(ii)$ there exists a solution
to $f(s,r,3)$.
\end{lemma}
\begin{proof}
Let $\Delta : [1, 6s - 5]$ $([1, s + 3r - 2 - \delta])
\longrightarrow \{1,2,3\} $ be given.  If one color occurs at most
 $s-1$ times, and if $A$ is the set of integers colored by the
 other two colors, then $|A| \geq 5s - 4$ $(|A| \geq 3r - 1 -
 \delta)$.  Hence, the proof is satisfied by Lemma 4.1.  Thus,
 each color occurs at least $s$ times.

Suppose that at most one color, call it the third color, denoted
by 3, contains at least $r$ elements.  Then, there are at least
 $6s - 5 - 2(\lfloor \frac{5s - 1}{3} \rfloor - 2) \geq 2 \lceil
 \frac{4s + 1 }{3} \rceil - 1 \geq s + r$ elements colored by the third
 color.  $($There are at least $s + 3r - 2- \delta - 2(r-1) \geq s
 + r - \delta$ elements colored by the third color.$)$  There can
 be at most $3s - 3$ integers in the interval $[first(3),
 last(3)]$, else condition $(i)$ holds.

Suppose that $r < \lfloor \frac{3s}{2} \rfloor - 1$.  Since the
first two colors each contain fewer than $r$ elements, the third
color consists of at least $6s - 5 - 2( \lfloor \frac{3s}{2}
\rfloor - 1) - 1 \geq 3s - 2$ elements.  Consider the set $B =
first_1^{r-1}(3) \cup last(3)$.  Since $diam(B) \geq 3s - 2$,
condition $(i)$ holds, and hence we can assume that $r \geq
\lfloor \frac{3s}{2} \rfloor - 1$.  There are at least $s + r - 1
\geq \lfloor \frac{5s}{2} \rfloor$ elements colored by the third
color.  Let $C = first_1^s(3)$.  Since there are at most $3s - 3$
elements colored by the third color else condition $(i)$ holds,
$diam(C) \leq \lceil \frac{3s}{2} \rceil - 2$.  Since $r \geq
\lfloor \frac{3s}{2} \rfloor - 1$ condition $(ii)$ holds by taking
the sets $C$ and $D = first_{s + 1}^{s + r - 1}(3) \cup last(3)$.
Thus, at least two of the three colors must contain at least $r$
elements.

Let w.l.o.g. $\Delta(1) = 1$.  Let the next new color be denoted
$2$ and the last new color be denoted $3$.  Suppose that the first
$s$ elements are monochromatic.  In this case, there can be at
most $r - 1$ more elements of each color else condition $(ii)$
holds.  Since $s + 3(\lfloor \frac{5s - 1}{3} \rfloor - 2) < 6s -
5$ $(s + 3(r-1) - 1 - \delta < s + 3r - 2 - \delta)$, the next
element in the set after $r-1$ elements of each color in any order
will obtain condition $(ii)$.  Hence, there exists some element of
$[1,s]$ not colored by color 1.  To complete the proof of the
lemma, we will consider four different cases based upon the
cardinality of colors 1 and 2. \newline \noindent \textbf{Case 1:}
Assume there are at most $r-1$ elements of color 1 and color 2.

Since there must be at least two colors with at least $r$
elements, we have a contradiction. \newline \noindent \textbf{Case
2:} Assume there are at least $r$ elements of color 1 and at most
$r-1$ elements of color $2$.

There must be at least $r$ elements of color $3$ because there
must be two or more colors with at least $r$ elements.  We will
now show that the interval $[1, 2s - 1]$ is only colored by colors
1 and 2.  Suppose that $first(3) \leq \lceil \frac{4s -
2}{3}\rceil - 1$.  Since color $3$ is colored by at least $r$
elements, $last(3) \leq \lceil \frac{13s - 2}{3}\rceil - 5$ else
condition $(i)$ holds via the set $E = first_1^{r-1}(3) \cup
last(3).$  Hence, the interval $[\lceil \frac{13s - 2}{3}\rceil -
4, 6s - 5]$ $([\lceil \frac{13s - 2}{3}\rceil - 4, s + 3r - 2 -
\delta]),$ consists of at least $\lceil \frac{5s}{3} \rceil - 1
\geq r$ $(3r - \lfloor \frac{10s}{3} \rfloor + 3 - \delta \geq r)$
elements colored by $2$, which is a contradiction.

We will repeat this argument two more times by shifting the
location of some elements colored by $2$ from the beginning of the
interval to the end of the interval.  Suppose that $\lceil
\frac{4s - 2}{3}\rceil - 1 < first(3) \leq \lceil \frac{5s -
2}{3}\rceil - 1$.  Hence, the interval $[\lceil \frac{14s -
2}{3}\rceil - 4, 6s -5]$ $([\lceil \frac{14s -2}{3} \rceil - 4, s
+ 3r - 2 - \delta])$ is colored entirely by $2$, which contains at
least $\lceil \frac{4s - 2}{3}\rceil - 1$ $(3r - 3s + 2 - \delta)$
elements.  There are at least $\lceil \frac{s - 2}{3}\rceil$
elements colored by $2$ in the interval $[1, \lceil \frac{4s -
2}{3}\rceil - 1]$, else we obtain a monochromatic set $F =
first_1^s(1, [1, \lceil \frac{4s - 2}{3}\rceil - 1])$ where
$diam(F) \leq \lceil \frac{4s - 2}{3}\rceil - 2 \leq r$, assuming
$r > \lceil \frac{4s - 2}{3}\rceil - 2$, thereby forcing condition
$(i)$.  Notice that if $r \leq \lceil \frac{4s - 2}{3}\rceil - 2$,
then there are already too many elements colored by $2$ in $[
\lceil \frac{14s - 2}{3}\rceil - 4, 6s - 5]$ $([ \lceil \frac{14s
- 2}{3}\rceil - 4, s + 3r - 2 - \delta])$.  Hence, condition
$(ii)$ holds by taking $E$ and $F$.


Next, suppose that $ \lceil \frac{5s - 2}{3}\rceil - 1 < first(3)
\leq 2s - 1.$  In this case, the interval $[1,  \lceil \frac{5s -
2}{3}\rceil - 1]$ requires $\lceil \frac{2s - 2}{3}\rceil $
elements colored by $2$ else condition $(ii)$ holds, assuming that
$ r \geq \lceil \frac{5s - 2}{3}\rceil - 1$.  Otherwise, if $r
\leq \lceil \frac{5s - 2}{3}\rceil - 1$, we know there are at
least $r - s + 1$ elements colored by $2$ in $[1,  \lceil \frac{5s
- 2}{3}\rceil - 1]$, else we can force condition $(ii)$ with $F' =
first_1^s(1) \cup last_1^r(3)$. Since the interval $[ \lceil
\frac{15s - 2}{3}\rceil - 4, 6s - 5]$ $([ \lceil \frac{15s -
2}{3}\rceil - 4, s + 3r - 2 - \delta])$ is colored entirely by
$2$, we require an additional $s - 1$ $(3r - \lfloor \frac{8s}{3}
\rfloor + 3 - \delta)$ elements for a total of $r - s + 1 + s - 1
= r$ $(3r - \lfloor \frac{10s - 2}{3} \rfloor + 2 - \delta \geq
r)$ elements colored by $2$, which is a contradiction.  Hence, the
interval $[1, 2s - 1]$ is colored entirely by $1$ and $2$.

By Theorem 1.1, there exists an $s$-element monochromatic set $G$
in $[1, 2s - 1]$ with $diam(G) \leq 2s - 2$.  (Thus, if $r \geq 2s
- 1$, condition $(ii)$ follows).  This implies that there exists
an interval in $[2s , 6s - 5]$ $([2s, s + 3r - 2- \delta])$ of
diameter at most $2s - 3$ consisting of all the elements colored
by $3$.  Since the interval $[3s - 2, 6s - 5]$ $([3s - 2, s + 3r -
2 - \delta])$ is colored only by $2$ and $3$, there must be at
least $s$ elements colored by $2$ in this interval.  Since there
are at least $r$ elements colored by $3$, $diam(G) \geq r$, else
condition $(ii)$ holds.  Hence, there are at least $r - s + 1$
elements colored by $2$ in the interval $[1, 2s - 1]$.  This
implies that there are at least $s + r - s + 1= r + 1$ elements
colored by 2, a contradiction, which completes this case. \newline
\noindent  \textbf{Case 3:} Assume there are at least $r$ elements
of color $1$ and color $2$.

Since $first(1) = 1, last(1) \leq 3s - 3$, else condition $(i)$
follows.  Similarly, since $first(2) \leq s$, then $last(2) \leq
4s - 4.$  Thus, the interval $[4s - 3, 6s - 5]$ $([4s - 3, s + 3r
- 2 - \delta])$ is colored entirely by color 3.  It follows that
there is a monochromatic $s$-element set $H \subset [1, 2s - 1]$
with $diam(H) \leq 2s - 2$, which along with $H' =
first_1^{r-1}(3) \cup last(3)$, yields condition $(ii)$.
\newline \noindent \textbf{Case 4:} Assume there are at most $r -
1$ elements of color $1$ and at least $r$ elements of color 2.

Since $first(2) \leq s$, then $last(2) \leq 4s - 4$, else
condition $(i)$ follows.  Hence the interval, $[4s - 3, 6s - 5]$
$([4s - 3, s + 3r - 2 - \delta])$ is entirely colored by 1 and 3.
Since at least two of the colors must contain at least $r$
elements, there must be at least $r$ elements colored by color 3.
Using the same argument as in Case 2 except that we condition on
the location of the elements colored by 1, the interval $[1, 2s
-1]$ is colored entirely by colors 1 and 2.  Hence, there exists a
monochromatic $s$-element set $I = first_1^s(1, [1, 2s - 1])$ or
$first_1^s(2, [1, 2s -1])$, such that $diam(I) \leq 2s - 2$.  If
there are $2s - 1$ or more elements colored by 3, then condition
$(ii)$ follows by taking $I$ and $J = first_1^{r-1}(3) \cup
last(3)$.  Thus, there are at most $2s - 2$ elements colored by 3.
In addition, it follows that there exists an interval of diameter
at most $2s - 3$ consisting of all the elements colored by $3$.
(Notice that condition $(ii)$ follows using $I$ and $J$ if $r \geq
2s - 1.)$

Suppose that $2s - 1 < first(3) \leq 3s - 2$.  Since there are at
least $r$ elements colored by $3$, there must be at least $r - s +
1$ elements colored by $1$ in $[1, 2s - 1]$ else condition $(ii)$
holds by taking $K = first_1^s(1, [1, 2s - 1])$ or $first_1^s(2,
[1, 2s -1])$ and $L = first_1^{r-1}(3) \cup last(3)$.  Since
$first(3) \leq 3s - 2, last(3) < 5s - 4$.  Hence, the interval
$[5s - 4, 6s - 5]$ $([5s - 4, s + 3r - 2 - \delta])$, consisting
of at least $s - 1$ elements is colored entirely by color $1$.  A
total of $(s - 1) + (r - s + 1) = r$ elements are colored by $1$,
a contradiction.  Hence, the interval $[1, 3s - 2]$ is colored
entirely by $1$ and $2$.

Next, suppose $r \leq \lfloor \frac{3s}{2} \rfloor - 2$.  If there
are at least $ \lfloor \frac{5s}{2} \rfloor - 1$ elements colored
by $2$, then we can construct $M = first_1^s(2)$ with $diam(M)
\leq \lceil \frac{3s}{2} \rceil -2$.  Condition $(ii)$ holds with
$M$ and $N = first_{s+1}^{s + r - 1}(2) \cup last(2)$, where
$diam(N) \geq \lfloor \frac{3s}{2} \rfloor - 2$.  Thus, there are
at most $\lfloor \frac{5s}{2} \rfloor - 2$ elements colored by 2,
and there are at least $\lceil \frac{3s}{2} \rceil   - 1$ elements
colored by 1.  (Recall $r < 2s - 1$ and $|\Delta^{-1}(3)| \leq 2s
- 2.$)  But this means that there are at least $r$ elements
colored by 1, a contradiction.  Therefore $r > \lfloor
\frac{3s}{2} \rfloor - 2$.

Suppose that $3s - 2 < first(3) \leq \lfloor \frac{10s}{3} \rfloor
- 3$.  Let $T = first_1^{r-1}(3) \cup last(3)$, where $diam(T) =
t$, and let $\alpha = 2s - 2 - t$.  Then the interval, $[\lfloor
\frac{16s}{3} \rfloor - 5 - \alpha, 6s - 5]$ $([\lfloor
\frac{16s}{3} \rfloor - 5 - \alpha, s + 3r - 2 - \delta]),$ is
colored entirely by $1$.  In this interval, there are at least $
\lceil \frac{2s}{3} \rceil + \alpha + 1 $ $(\lceil \frac{2s}{3}
\rceil + \alpha + 1 + 3(r - \lfloor \frac{5s - 1}{3} \rfloor) - 1
- \delta)$ elements colored by $1$.  Thus, there are at most $s -
2 - \alpha$ elements left that can be colored by $1$.  Consider
the interval $[2, 2s - \alpha]$.  Since $\Delta(1) = 1$, there are
$s - 3 - \alpha$ elements of color 1 remaining before a
contradiction occurs.  We must also avoid creating an $s$-element
monochromatic set of diameter $t$ or less else condition $(ii)$
holds.  Thus, there must be at least $t - s + 1  = s - 1 - \alpha$
elements colored by $1$ in the interval $[2, 2s - \alpha]$.
However, this is more than $s - 3 - \alpha$ elements, and so we
have obtained a contradiction.  Thus, $first(3) > \lfloor
\frac{10s}{3} \rfloor - 5 $.

Next, suppose that $\lfloor \frac{10s}{3} \rfloor - 3 < first(3)
\leq \lfloor \frac{11s}{3} \rfloor - 3$.  Then the interval
$[\lfloor \frac{17s}{3} \rfloor - 5 - \alpha, 6s - 5]$  $( [
\lfloor \frac{17s}{3} \rfloor - 5 - \alpha, s + 3r - 2 - \delta])$
is colored entirely by $1$ and the interval $[1, \lfloor
\frac{10s}{3} \rfloor - 3]$ is colored entirely by $1$ and $2$. In
this interval, there are at least $\lfloor \frac{s}{3} \rfloor +
\alpha + 1$ $(\lceil \frac{s}{3} \rceil + \alpha + 1 + 3(r -
\lfloor \frac{5s - 1}{3} \rfloor) - 1 - \delta)$ elements colored
by 1. Thus, there are at most $\lfloor \frac{4s}{3} \rfloor - 2 -
\alpha$ elements left that can be colored by 1 before we have a
contradiction.  The diameter of the elements colored by color $2$
is at most $3s - 3$, so there are $\lfloor \frac{s}{3} \rfloor$
elements that are forced to be colored by 1 outside the interval
$[first(2), last(2)]$ in $[1, \lfloor \frac{10s}{3} \rfloor - 3
]$.  Subtracting those elements colored by $1$ leaves $s - 2
-\alpha$ integers that are allowed to be colored by 1.  Inside
$[first(2), last(2)]$, we must avoid creating an $s$-element
monochromatic set of diameter $t$ or less, else condition $(ii)$
holds.  Thus, there must be at least $t - s + 2 = s + \alpha$
elements  colored by $1$ in the interval $[first(2), first(2) + 2s
- 2 - \alpha]$. But there aren't enough elements colored by $1$ to
satisfy these constraints.  Thus, $ first(3) > \lfloor
\frac{11s}{3} \rfloor - 3$.

Next suppose that $\lfloor \frac{11s}{3} \rfloor - 3 < first(3)
\leq 4s - 3$.  Hence the interval $[6s - 5 - \alpha, 6s - 5]$
$([6s - 5 - \alpha, s + 3r - 2 - \delta])$ is colored entirely by
$1$.  In this interval, there are at least $\alpha + 1$ $(\alpha +
1 + 3(r - \lfloor \frac{5s-1}{3} \rfloor) - 1 - \alpha)$ elements
colored by 1.  thus, there are at most $\lfloor \frac{5s}{3}
\rfloor - 2 - \alpha$ elements that are colored by 1.  The
diameter of the elements colored by color 2 is at most $3s - 3$,
so there are $\lfloor \frac{2s}{3} \rfloor$ elements that must be
colored by $1$ outside this interval.  Subtracting those elements
of 1 leaves $s - 2- \alpha$ integers colored by 1 before we obtain
a contradiction. We must avoid creating an $s$-element
monochromatic set of diameter $t$ or less in the interval
$[first(2), first(2) + 2s - 2- \alpha]$, else condition $(ii)$
holds. Thus, there must be at least $t -s  + 2 = s - \alpha$
elements colored by $1$ in this interval.  But we only have $s - 2
- \alpha$ elements colored by 1 left before we obtain a
contradiction.  Thus, $first(3) > 4s - 3$.

Finally, suppose that $first(3) > 4s - 3$.  Then the intervals
$[4s - 3, first(3)  - 1]$ and $[last(3) + 1,6s - 5] $ $([last(3) +
1, s + 3r - 2 - \delta])$ are colored entirely by $1$.  In the
interval $[4s - 3, 6s - 5]$ $([4s - 3, s + 3r - 2 - \delta])$
least, there are at least $\alpha$ $(\alpha + 3(r - \lfloor
\frac{5s-1}{3} \rfloor) - 1 - \delta)$ elements colored by $1$.
Thus, there are at most $\lfloor \frac{5s}{3} \rfloor - 2 -
\alpha$ elements left that may be colored by 1.  The diameter of
the elements colored by color 2 is at most $3s - 3$, so there are
(at least) $s - 1$ elements that must be colored by 1 in $[1, 4s -
3]$ outside $[first(2), last(2)]$ because the interval $[1, 4s -
3]$ is colored entirely by $1$ and $2$.  Subtracting those
elements of $1$ leaves at most $\lfloor \frac{2s}{3} \rfloor - 1 -
\alpha$ integers colored by 1.  We must avoid creating an
$s$-element monochromatic set of length $t$, else condition $(ii)$
holds. Thus, there must be at least $t - s - 2 = s - \alpha$
elements colored by $1$ in the interval, $[first(2), first(2) + 2s
- 2 - \alpha]$.  However, there aren't enough elements colored by
$1$ to satisfy these constraints without obtaining a contradiction
to our assumptions.  Thus, our proof of the lemma is complete.
\qed
\end{proof}

\begin{theorem}
Let $r$ and $s$ be positive integers where $r \geq s \geq 3$.  If
$r \leq \lfloor \frac{5s-1}{3} \rfloor - 1$ then $f(s, r, 3) = 9s
- 7.$
\end{theorem}
\begin{proof}
The coloring $\Delta : [1, 9s - 8] \longrightarrow \{1,2,3\}$
given by the string:
\begin{align*}
12^{s-1}1^{s-2}3^{s-1}12^{s-1}1^{s-1}2^{s-1}1^{s-1}3^{2s-2}
\end{align*}
implies that $f(s,r,3) > 9s - 8$.  Next we show that $f(s,r,3)
\leq 9s -7.$  Let $\Delta : [1, 9s-7] \longrightarrow \{1,2,3\}$
be an arbitrary coloring.  By Theorem 1.1, the interval $[1, 3s -
2]$ contains an $s$-element monochromatic subset $S$ such that
$diam(S) \leq 3s - 3$.  Shifting the interval $[3s - 1, 9s - 7]$
to the interval $[1, 6s - 5]$ and applying Lemma 4.2 completes the
proof for $r \geq s \geq 3$. \qed
\end{proof}

\begin{theorem}
Let $r$ and $s$ be positive integers where $r \geq s \geq 3$.  If
$ \lfloor \frac{5s-1}{3} \rfloor - 1 < r \leq 2s - 2$ then $f(s,
r, 3) = 4s + 3r - 4.$
\end{theorem}
\begin{proof}
The coloring $\Delta : [1, 4s + 3r - 5] \longrightarrow \{1,2,3\}$
given by the string:
\begin{align*}
12^{s-1}1^{s-2}3^{s-1}12^{s + r - 1}1^{r-1}3^{r-1}
\end{align*}
implies that $f(s,r,3) > 4s + 3r - 5$.  Next we show that
$f(s,r,3) \leq 4s + 3r - 4.$  Let $\Delta : [1, 4s + 3r - 4]
\longrightarrow \{1,2,3\}$ be an arbitrary coloring.  By Theorem
1.1, the interval $[1, 3s - 2]$ contains an $s$-element
monochromatic subset $S$ such that $diam(S) \leq 3s - 3$. Shifting
the interval $[3s - 1, 4s + 3r - 4]$ to the interval $[1, s + 3r -
2]$ and applying Lemma 4.2 completes the proof for $r \geq s \geq
3$. \qed
\end{proof}

\begin{theorem}
Let $r$ and $s$ be positive integers where $r \geq s \geq 3$.  If
$ 2s - 2 < r \leq 3s - 3$ then $f(s, r, 3) = 6s + 2r - 6.$
\end{theorem}
\begin{proof}
The coloring $\Delta : [1, 6s + 2r - 7] \longrightarrow \{1,2,3\}$
given by the string:
\begin{align*}
12^{s-1}1^{s-2}3^{s-1}12^{3s - 3}1^{r-1}3^{r-1}
\end{align*}
implies that $f(s,r,3) > 6s + 2r - 7$.  Next we show that
$f(s,r,3) \leq 6s + 2r - 6.$  Let $\Delta : [1, 6s + 2r - 6]
\longrightarrow \{1,2,3\}$ be an arbitrary coloring.  By Theorem
1.1, the interval $[1, 3s - 2]$ contains an $s$-element
monochromatic subset $S$ such that $diam(S) \leq 3s - 3$. Shifting
the interval $[3s - 1, 6s + 2r - 6]$ to the interval $[1, 4s + 2r
- 4]$ and applying Lemma 4.2 completes the proof for $r \geq s
\geq 3$. \qed
\end{proof}

\begin{theorem}
Let $r$ and $s$ be positive integers where $r \geq s \geq 2$.  If
$  r > 3s - 3$ then $f(s, r, 3) = 3s + 3r - 4.$
\end{theorem}
\begin{proof}
The coloring $\Delta : [1, 3s + 3r - 5] \longrightarrow \{1,2,3\}$
given by the string:
\begin{align*}
12^{s-1}1^{s-2}3^{s-1}12^{r-1}1^{r-1}3^{r-1}
\end{align*}
implies that $f(s,r,3) > 3s + 3r - 5$.  By the pigeonhole
principle, it follows that the interval $[1, 3s - 2]$ contains an
$s$-element monochromatic subset $S_1$ with $diam(S_1) \leq 3s -
3$, and that the interval $[3s - 1, 3s + 3r - 4]$ contains an
$r$-element monochromatic subset $S_2$ with $diam(S_2) \geq 3s -
3$.  Hence, $S_1$ and $S_2$ complete the proof. \qed
\end{proof}
Now we will determine $f(s, r, \{\infty\} \cup \mathbb{Z})$.

\begin{lemma}
Let $r$ and $s$ be positive integers where $r \geq s \geq 3$. If
$\Delta : [1, 6s - 5 + \frac{s}{(r,s)} - 1] \longrightarrow \{
\infty \} \cup \mathbb{Z}$ $ (\Delta : [1, s + 2r - 4 + \min(3s -
3, s + r - 1 + \frac{s}{(r,s)} -1)] \longrightarrow \{ \infty \}
\cup \mathbb{Z})$ is a coloring, then either:
\newline \noindent $(i)$ there exists an $r$-element subset $S
\subseteq [1, 6s - 5 + \frac{s}{(r,s)} - 1]$ $( S \subseteq [1, 2r
- 1 + \min(3s - 3, s + r - 1 + \frac{s}{(r,s)} -1)])$ such that

$(a)$ $S$ is either $\infty$-monochromatic or zero-sum mod $r$ and

$(b)$ $diam(S) \geq 3s - 3$ or \newline $(ii)$ there exists a
solution to $f(s, r, \{ \infty \} \cup \mathbb{Z})$.
\end{lemma}
\begin{proof}
We partition the proof into three cases, each considering a
different cardinality of the sets $\Delta^{-1}(\mathbb{Z})$ and
$\Delta^{-1}(\infty)$.
\newline
\noindent \textbf{Case 1:} Suppose $|\Delta^{-1}(\mathbb{Z})| \leq
3s - 3$.

This implies that $|\Delta^{-1}(\infty)| \geq 3s - 2$. There
exists an $r$-element subset $S \subseteq \Delta^{-1}(\infty)$
such that $S$ is $\infty$-monochromatic and $diam(S) \geq 3s - 3$.
Condition $(i)$ of the lemma follows.
\newline
\noindent \textbf{Case 2:} Suppose that $|\Delta^{-1}(\infty)|
\leq r - 1$. \newline \noindent \textbf{Case 2.1:} First, consider
the case when $r \leq \lfloor \frac{3s}{2} \rfloor - 1 $ or $2s -
2 \leq r \leq 3s - 3$.

In this case we do not have to apply Theorem 2.3.  Let $I_1 =
first_1^r(\mathbb{Z})$ and let $I_3 = last_1^r(\mathbb{Z})$.  Also
let $a = last(\mathbb{Z}, I_1), b = first(\mathbb{Z}, I_3), S_1 =
(I_1 \setminus \{a\}) \cup I_3$ and $S_2 = I_1 \cup (I_3 \setminus
\{b\})$.  Since $|S_1| = |S_2| = 2r - 1$, by Theorem 1.1, there
exist $r$-element zero-sum mod $r$ sets $X, Y$ in $S_1, S_2$,
respectively.  If $X \cap (I_1 \setminus \{a\}) \neq \emptyset$
and $ X \cap I_3 \neq \emptyset$, then $diam(X) \geq 3s - 3$, and
assertion $(i)$ of the lemma is satisfied.  Else, $X = I_3$ since
$|I_1 \setminus \{a\}| \leq r$.  Repeating the above argument for
$Y$, either $(i)$ follows or $Y = I_1$.  Hence, $I_1$ and $I_3$
are $r$-element zero-sum mod $r$ sets.

Consider the coloring of the sets $I_1$ and $I_3$.  If
$\Delta(I_1) \cap \Delta(I_3) \neq \emptyset$, then replace an
element $u \in I_1$ with an element $v \in I_3$, where $\Delta(u)
= \Delta(v)$ mod $r$ to construct an $r$-element zero-sum mod $r$
set $S$ with $diam(S) \geq 3s - 3$.  Hence, assume that
$\Delta(I_1) \cap \Delta(I_3) = \emptyset$.  Suppose that at least
one of the sets $I_1$ and $I_3$ is not monochromatic.  Without
loss of generality, suppose that $I_1$ is colored by at least two
residues mod $r$.  Let $M = first_1^{r-1}(\mathbb{Z}, I_1)$, and
let $M' = last_1^{r-1}(\mathbb{Z}, I_1)$.  Either $M$ or $M'$, say
$M$, is colored by at least two residues.  Let $N =
first_1^{r-1}(\mathbb{Z}, I_3)$.  Consider the set $A = M \cup N$.
It is colored by at least three residues mod $r$ since $M$ is
colored by at least two residues and $N$ is colored by at least
one residue.  By Theorem 2.1, there exists an $r$-element zero-sum
mod $r$ set $A' \subset A$ with $diam(A') \geq 3s - 3$, whence
condition $(i)$ follows.  Thus, both $I_1$ and $I_3$ are
monochromatic mod $r$.

Let $J_1 = I_1 \setminus \{a\}, J_2 = \{first(\mathbb{Z}, [first(
\mathbb{Z}) + 3s - 2, last(\mathbb{Z})])\}$, and $J_3 =
last_1^{r-1}(I_3)$.  Within $J_1 \cup J_2 \cup J_3$, there are $2r
- 1$ elements, and so there exists an $r$-element zero-sum mod $r$
set $B$.  If $B$ is contained in $B_1 = J_1 \cup J_2$, $B_2= J_1
\cup J_3$ or $B_3 = J_1 \cup J_2 \cup J_3$, with the additional
condition that each of the named component sets in each union must
be nonempty, (i.e. for $B_1$, both $J_1$ and $J_2$ must be
nonempty) then condition $(i)$ follows. If not, then $B$ is
contained in $J_2 \cup J_3$. If $r \geq s + \frac{s}{(r,s)} - 1$,
then in $I_1$, there exists an $s$-element zero-sum mod $s$ set
$C$.  If $\alpha$ counts the number of elements colored by
$\infty$ in $[first(B), last(B)]$, then $diam(B) \leq s +
\frac{s}{(r,s)} - 2 + \alpha$. If $r \geq s + \frac{s}{(r,s)} -
1$, then condition $(ii)$ follows by taking $B$ and $C$, where
$diam(C) \geq \lceil \frac{3s}{2} \rceil - 3 + \frac{s}{(r,s)} +
\alpha$.

$($If $r < s + \frac{s}{(r,s)} - 1$, then first suppose that there
are only two residues colored by $\mathbb{Z}$.  There cannot be
only one residue colored by $\mathbb{Z}$ since $\Delta(I_1) \cap
\Delta(I_3) = \emptyset$.  Let 1 be the residue mod $r$ used to
color $I_1$ and let $2$ be the residue mod $r$ used to color
$I_3$. This implies that both residue 1 and residue 2 contain at
least $r$ elements mod $r$.  Notice that $last(1) < J_2$ and
$first(2) > first_{3s - r- 1}(\mathbb{Z})$ else condition $(i)$
holds.  This implies that the interval $[1, first_{3s - r-
1}(\mathbb{Z})]$, is colored entirely by $\infty$ or $1$ mod $r$.
Also, the interval $[J_2, last(\mathbb{Z})]$ is colored entirely
by residue $2$ mod $r$.  If $\alpha$ is the number of elements
colored by $\infty$ in $[first(\mathbb{Z}), first_{3s - r-
1}(\mathbb{Z})]$, then there exists an $s$-element zero sum mod
$s$ set $D$ where $diam(D) \leq s + \frac{s}{(r,s)} - 1 + \alpha
\leq 3s - r - 1 + \alpha$. Since there are at least $3s - r - 1 +
\alpha$ elements colored by residue $2$ mod $r$, condition $(ii)$
holds by taking $D$ and $first(2) \cup last_1^{r-1}(2)$.)

\noindent \textbf{Case 2.2:} Next, consider the case when $\lfloor
\frac{3s-1}{2} \rfloor - 1 \leq r \leq 2s - 2$.

In this case, we must apply Theorem 2.3.  Let $a =
first(\mathbb{Z}), b = first_2(\mathbb{Z}), y =
last_2(\mathbb{Z})$ and $z = last(\mathbb{Z})$.  Since
$|\Delta^{-1}(\infty)| \leq r - 1$, then
$|\Delta^{-1}(\mathbb{Z})| \geq 6s - r - 5 + \frac{s}{(r,s)} - 1$
$(|\Delta^{-1}(\mathbb{Z})| \geq r + \omega$, where $\omega =
\min(3s - 3, s + r - 1 + \frac{s}{(r,s)}- 1)$). Consider the
interval $S = [a, z]$. Let $c$ be the residue of $\mathbb{Z}$ that
colors the most number of integers mod $r$, and let $d$ be the
residue that colors the second most number of integers mod $r$.

Suppose there are only two residues that are colored by
$\mathbb{Z}$ mod $r$.  If residue $d$ had fewer than $r$ elements
mod $r$, then residue $c$ would consist of at least $s + \lfloor
\frac{5s}{3} \rfloor +  \frac{s}{(r,s)} - 2$ $(1 + \omega)$
elements and condition $(ii)$ would follow.  Thus, both residues
must contain at least $r$ elements mod $r$.  This means that the
diameters of the set of elements containing each residue is at
most $3s -3$. Let $1$ be the residue that appears first and let
$2$ be the residue that appears second, making $first(1) <
first(2)$.  The set $A = first_1^{s + \frac{s}{(r,s)} - 1}(1)$
contains an $s$-element zero-sum mod $s$ set, $A'$.  Recall that
there are at least $s + 2r - 2 + \frac{s}{(r,s)}$ elements colored
by 1 or 2.  Notice that $first(2)$ must occur after $first_{s +
\frac{s}{(r,s)} - 1}(1)$, else condition $(i)$ follows because
there would be more than 3s - 3 elements occurring after both
residues mod $r$ have already appeared.  Thus $first(2)
> first_{s + \frac{s}{(r,s)} - 1}(1)$.  But then condition
$(ii)$ follows by taking the sets $A'$ and $B = first(2) \cup
last_1^{r-1}(2)$. To see that $diam(A) \leq diam(B)$, if $\alpha$
is the number of elements colored by $\infty$ within the interval
that contains $A'$, there must be at least $s + 2r - 1 +
\frac{s}{(r,s)} - 1 - (3s - 3 - \alpha) = 2r - 2s + 2 +
\frac{s}{(r,s)} + \alpha \geq s + \frac{s}{(r,s)} - 1$. Hence,
there must exist some integer that is colored by a third residue
mod $r$.\newline \noindent \textbf{Case 2.2.1} Suppose there are
at least $r + 1$ elements of $\mathbb{Z}$ not colored by $c$ mod
$r$.

We must satisfy the conditions of Theorem 2.3 so we can produce an
$r$-element zero-sum mod $r$ set $E$ with $diam(E) \geq 3s - 3$ if
conclusion $(i)$ of Theorem 2.3 holds.  If there are at least two
elements not colored by $c$ or $d$, we satisfy the conditions of
Theorem $2.3$, by 1) creating an $(r-2)$-set partition of
$\Delta(S)$ with $r - 3$ sets of cardinality $2$ and one set of
cardinality $3$, and 2) fixing two elements that are not part of
the set partition, in this case $a$ and $z$.  This holds unless
the only two elements not colored by $c$ or $d$ are located at $a$
and $z$, whence we cannot force a set with diameter $z - a$.  In
this case, fix $b$ and $z$, and the diameter of this set described
by Theorem 2.3 is still at least $3s - 3$.  Hence, we can apply
Theorem 2.3 to this case if there are at least two elements of
$\mathbb{Z}$ not colored by $c$ or $d$ mod $r$.

Suppose that there is only one element of $\mathbb{Z}$ colored by
something other than $c$ or $d$ mod $r$.  In this case, fix $a$
and $z$ for the set partition, except if either $a$ or $z$ is the
third color.  If $a$ is the third color, fix $b$.  Likewise, if
$z$ is the third color, $y$.  In any of these three instances, the
diameter of the $r$-element zero-sum mod $r$ set that may be
satisfied by Theorem 2.3 is still greater than $3s- 3$.  Hence we
can apply Theorem 2.3 to this case, and will do this after Case
2.2.2. \newline \noindent \textbf{Case 2.2.2:} Suppose that all
but $r$ elements of $\mathbb{Z}$ are colored by $c$ mod $r$.
Notice that if all but $r-1$ elements of $\mathbb{Z}$ are colored
by $c$, color $c$ would consist of at least $s + \lfloor
\frac{5s}{3} \rfloor + \frac{s}{(r,s)} - 2$ $(1 + \omega)$
elements, and condition $(ii)$ would follow.  Therefore, there
exists a $c \in \mathbb{Z}$ such that $|\Delta^{-1}(c) \cup S| =
6s - 2r - 5 + \frac{s}{(r,s)} \geq s + \frac{s}{(r,s)} + \lfloor
\frac{5s}{3} \rfloor - 3$ $(|\Delta^{-1}(c) \cup S| = \omega)$ mod
$r$.  Thus, if there are at least $\lceil \frac{s}{3} \rceil -
\frac{s}{(r,s)} + 1$ elements not colored by $c$ in $[first(c),
last(c)]$, condition $(i)$ follows since we have constructed an
$r$-element zero-sum mod $r$ set $E$ with $diam(E) \geq 3s - 3$.
Let $F = first_1^{s + \frac{s}{(r,s)} -1}(c)$. There exists an
$s$-element zero-sum mod $s$ set $F'$ such that $diam(F') \leq
\lceil \frac{4s}{3} \rceil \leq r$.  Take the other $r-1$ elements
colored by $c$ and combine them with $r-1$ elements colored by
$\mathbb{Z}$ having at least two different colors than $c$ mod
$r$.  These differently colored residues exist since we have
already considered the case when there are only two residues of
$\mathbb{Z}$ mod $r$.  Hence, in this set, there exists an
$r$-element zero-sum mod $r$ set $G$. If all the elements of $G$
are after $F'$, then condition $(ii)$ holds.  If there exist
elements of $G$ with $F'$ or before $F'$, we can construct an
$(r-2)$-set partition in the following manner: Fix $a$ and $z$. If
either of these elements is the only element that is colored by
the third color, then replace it with $b$ or $y$, respectively.
After fixing the above elements, the remaining elements form a
$(2r - 3)$-element $(r-2)$-set partition.  We are left with the
case when $r < 4$.  This only occurs when both $r$ and $s$ are 3
and this is considered in \cite{bial}, completing Case 2.2.2.

Therefore, in both Case 2.2.1 and Case 2.2.2, there exists an
$(r-2)$-set partition of $\Delta(S)$ with $m-3$ sets of
cardinality $2$ and one set of cardinality $3$ as well as two
fixed elements.  We apply Theorem 2.3 to observe that 1) either
there is an $r$-element zero-sum mod $r$ set $H$ that contains the
two fixed elements and $r-2$ other elements of the set partition,
one from each set, or 2) all but at most $r-2$ of the elements
colored by $\mathbb{Z}$ are colored by elements from the same
coset ($a\mathbb{Z}_m + \alpha)$ of $\mathbb{Z}_m$.  If the first
conclusion occurs, $diam(H) \geq 3s - 3$, whence condition $(i)$
follows.  In the latter case, Theorem 2.2 implies that any subset
of cardinality $(r + \frac{r}{a} - 1 + a - 2) \leq \lceil
\frac{3s}{2} \rceil - 1$ must contain an $r$-element zero-sum mod
$r$ set.  Hence, consider the set $J = first_1^{\lceil
\frac{3r}{4} \rceil}(\mathbb{Z}) \cup last_1^{ \lceil \frac{3r}{4}
\rceil}(\mathbb{Z})$.  Within $J$, there exists an $r$-element
zero-sum mod $r$ set $K$ with $diam(K) \geq 3s - 3$.  Condition
$(i)$ holds, and the proof for this case is complete. \newline
\noindent \textbf{Case 3:} $|\Delta^{-1}(\mathbb{Z})| \geq 2s - 1$
and $|\Delta^{-1}(\infty)| \geq r$.
\newline \noindent \textbf{Case 3.1:} Suppose there are only one or
two residues colored by $\mathbb{Z}$ mod $r$.

Clearly, if there is only one residue colored by $\mathbb{Z}$ mod
$r$, we are finished.  Now, suppose that there are precisely two
residues of $\mathbb{Z}$ mod $r$.  Without loss of generality, let
color $1$ be the first residue colored by $\mathbb{Z}$ mod $r$ and
let $2$ be the second color mod $r$.

\noindent \textbf{Case 3.1.1:} Suppose that $first(\infty) <
first(1)$.  Suppose that $first(1)> s$.  Then condition $(ii)$
holds by taking $first_1^s(\infty)$ along with an $r$-element
zero-sum mod $r$ set from $last_1^{2r-1}(\mathbb{Z})$.  Therefore,
let $first(1) \leq s$.  Let $A_1 = first_1^{r-1}(\mathbb{Z}), A_2
= last_1^{r-1}(\mathbb{Z})$ and $A_3 = first(\mathbb{Z},
[first(\mathbb{Z}) + 3s - 3, last(\mathbb{Z})])$.  Notice that if
$A_3$ is contained in $A_2$, then $first(\mathbb{Z}) > s$, which
is a contradiction.  By Theorem 1.1, there exists an $r$-element
zero-sum mod $r$ set $A$.  If $A$ is contained using only and
precisely $C_1$ and $C_2$, $C_1$ and $C_3$ or elements from all
three sets, then condition $(i)$ holds.

Otherwise, $A$ contains elements from only $C_2$ and $C_3$.
Suppose that $C_1$ or $C_3$ is not monochromatic mod $r$.  If this
were so, we could interchange an element of $C_1$ with an element
in $C_2$ or $C_3$, and condition $(i)$ would follow.  Therefore,
$C_1$ and $C_3$ are monochromatic mod $r$.  Let $\alpha =
first(1)$.  Then $first(2) \leq 3s - 3 + \alpha$.  So
$|\Delta^{-1}| \geq 6s - 5 + \frac{s}{(r,s)} - 1 - (3s - 4 +
\alpha) \geq 3s - 1 - \alpha + \frac{s}{(r,s)} - 1$.  However,
notice that in the interval $[1, 2s + \frac{s}{(r,s)} - 2]$, there
exists either an $s$-element $\infty$-monochromatic set or an
$s$-element zero-sum mod $s$ set of diameter at most $2s - 2 +
\frac{s}{(r,s)} - 1$.  Call this set $B$.  The maximum value of
$\alpha$ is $s$, so $|\Delta^{-1}(2)| \geq 2s -1 + \frac{s}{(r,s)}
- 1$, and so the minimum diameter of $C$, an $r$-element zero-sum
mod $r$ set, colored entirely by $2$ is $2s - 2 + \frac{s}{(r,s)}
- 1$.  Therefore, $B$ and $C$ satisfy condition $(ii)$, completing
this subcase.

\noindent \textbf{Case 3.1.2:} Suppose that $first(1) <
first(\infty)$.  Now, suppose that $first(\infty) > 2s - 1$.  Then
there are at most $2s - 1$ colored by $\infty$, and $last(\infty)
- first(\infty) \leq 2s - 2$.  By the pigeonhole principle, one of
the two residues of $\mathbb{Z}$ mod $r$ must contain at least $r$
elements.  Call this residue $R_1$ and call the other residue
$R_2$, which may or may not contain $r$ elements.

\noindent Step 1:  If there is only one residue mod $r$ in $[1, s
+ \frac{s}{(r,s)} - 1]$, then condition $(i)$ or $(ii)$ holds
using only residues of $\mathbb{Z}$ when $ r < s + \frac{s}{(r,s)}
- 1$.  Otherwise, condition $(ii)$ holds by taking $D$, the
$s$-element zero-sum mod $s$ set in $[1, \frac{s}{(r,s)} - 1]$ and
$E$, the $r$-element $\infty$-monochromatic set in
$first_1^r(\infty)$.  Therefore, residues $R_1$ and $R_2$ mod $r$
occur before $s + \frac{s}{(r,s)} - 1$.  We can then conclude that
$last(R_1) \leq 4s - 4 +\frac{s}{(r,s)} $.  Therefore, the
interval $[4s - 3 + \frac{s}{(r,s)}, 6s - 6 + \frac{s}{(r,s)} ]$
$([4s - 3 + \frac{s}{(r,s)}, \max(last(\mathbb{Z}), last(\infty))
])$ is colored entirely by $\infty$ and $R_2$.

\noindent Step 2:  Suppose that $2s \leq first(\infty) \leq 3s
-3$. There are at most $2s - 2$ elements colored by $\infty$, else
condition $(ii)$ would hold.  Also, $last(\infty) \leq 5s - 6$.
Let $\gamma = 2s - 3 - (last(\infty) - first(\infty))$.  Notice
that there are at least $s + \frac{s}{(r,s)} + \gamma$ elements
colored by $R_2$ after $last(\infty)$.  Therefore there can be at
most $r - s - \frac{s}{(r,s)} $ elements colored by $R_2$ in the
rest of the interval else condition $(i)$ holds.  However, we know
that there are at least $2s - 1 - \gamma - (s + \frac{s}{(r,s)} -
1) = s - \frac{s}{(r,s)} - \gamma$ elements colored by $R_2$ in
$[1, 2s - 1]$, else condition $(ii)$ holds using an $s$-element
zero-sum mod $s$ set contained in $[1, 2s - 1]$ using elements
only from $R_1$ and an $r$-element $\infty$-monochromatic set. But
this means there are at least $s + \gamma + \frac{s}{(r,s)} + s -
\gamma - \frac{s}{(r,s)} = 2$ elements colored by $R_2$ and
condition $(ii)$ holds.  Notice that if $r > 2s$, then there are
more than $2s - 2$ elements colored by $\infty$ and condition
$(ii)$ holds using an $s$-element zero-sum mod $s$ set in the
interval $[1, 2s - 1]$ and the $r$ element $\infty$-monochromatic
set that is of diameter at least $2s$.

Step 3: Suppose that $3s - 2 \leq first(\infty)$.  Again, there
are at most $2s - 2$ elements colored by $\infty$, else condition
$(ii)$ would hold. Let $\gamma = 2s - 3 - (last(\infty) -
first(\infty))$.  Outside of the intervals $[first(R_1),
last(R_1)]$ and $[first(\infty), last(\infty)]$, there are at
least $s + \frac{s}{(r,s)} - 1 + \gamma$ elements colored by
$R_2$.  Inside the interval $[first(R_1), last(R_1)]$, there are
an additional $s - \frac{s}{(r,s)} - \gamma + 1$, else condition
$(ii)$ holds by applying an analogous argument from the previous
step.  This means there are at least $2s$ elements colored by
$R_2$ and condition $(i)$ holds if $r < 2s$ and the last element
is not colored by $\infty$.  If $r > 2s$, then there are more than
$2s - 2$ elements colored by $\infty$, a contradiction.  If the
last element is colored by $\infty$, condition $(i)$ still holds
since $first(R_2) \leq s + \frac{s}{(r,s)} - 1$ and $last(R_2)
\geq 4s - 3 + \frac{s}{(r,s)} - 1$, meaning the diameter of this
$r$-element zero-sum mod $r$ set has diameter at least $3s -3$.

Step 4:  Suppose that $s + \frac{s}{(r,s)} - 1 < first(\infty)
\leq 2s- 1$.  Again, both $R_1$ and $R_2$ must be contained in the
interval $[1, s + \frac{s}{(r,s)} - 1]$.  Therefore, the interval
$[5s - 3, \max(last(\mathbb{Z}), last(\infty)]$, is colored
entirely by $R_2$.  Thus, if there are at least $r$ elements of
$R_2$, condition $(i)$ holds.  Therefore, there are fewer than $r$
elements colored by $R_2$.

Consider two subcases.  First, suppose that $last(\infty) -
first(\infty) \geq first(\infty) - 1$.  Let $\lambda =
first(\infty) - (s + \frac{s}{(r,s)} - 1 ).$  This is the least
number of elements colored by $R_2$ before $first(\infty)$.  In
addition, there are an additional $5s - 3 - (3s - 3) -
first(\infty) = 2s - first(\infty)$ elements colored by $R_2$
which are after $[first(\infty), last(\infty)]$ but before $5s -
3$.  Therefore there are at least $s - 1 + \frac{s}{(r,s)} - 1 +
first(\infty) - (s + \frac{s}{(r,s)} - 1 ) + 2s - first(\infty) =
2s$ elements colored by $R_2$, and condition $(i)$ follows.

Now, suppose that $last(\infty) - first(\infty) < first(\infty) -
1$.  This means that $last(\infty) - first(\infty) \leq 2s - 2$.
the interval $[4s - 3 + \frac{s}{(r,s)} - 1,
\max(last(\mathbb{Z}), last(\infty))]$ is colored entirely by
$R_2$, and condition $(i)$ holds using the residues of $R_2$ mod
$r$.

\noindent Step 5: Suppose that $2 \leq first(\infty) \leq
\frac{s}{(r,s)} - 1$.  Consider three cases.  In all cases,
suppose that there is at least one color with $r$ elements mod
$r$, and call this residue $R_3$.  Call the other residue, which
may or may not contain $r$ elements mod $r$, $R_4$.  If both $R_3$
and $R_4$ are present before $first(\infty)$, then the last $2s -
2$ elements are colored by $R_4$.  Then condition $(i)$ holds
using the elements in $R_4$.  (If $r \geq 2s - 1$, then the last
$4s - 3$ elements are colored by $R_4$ and condition $(i)$ holds.)

Now suppose that $R_3$ is the only residue colored before
$first(\infty)$.  Again, the last $2s  -2$ elements must be
colored by $R_4$.  Therefore, the interval $[1, 3s - 3 +
\frac{s}{(r,s)} - 1]$ is colored entirely by $R_3$ and $\infty$.
Suppose that $first(\infty) \leq s - 1$.  Then $last(\infty) \leq
4s - 4$.  Thus, the interval $[4s - 3, \max( last(\mathbb{Z}),
last(\infty))]$ is colored entirely by $R_4$. However, in the
interval $[1, 2s - 1 + \frac{s}{(r,s)} - 1]$, there either exists
an $s$-element zero-sum mod $s$ set or an $s$-element
$\infty$-monochromatic set.  This set has diameter at most $2s - 2
+ \frac{s}{(r,s)} - 1 $.  Similarly, in $[4s - 3, \max(
last(\mathbb{Z}), last(\infty))]$, there exists an $r$-element
zero-sum mod $r$ set with diameter at least $2s - 2 +
\frac{s}{(r,s)} - 1$.  Therefore, condition $(ii)$ holds.

So let $s \leq first(\infty) \leq s + \frac{s}{(r,s)} - 1$. Notice
that if $first(R_4) \leq 4s -4$, then condition $(ii)$ will hold
in the following manner:  There exists an $r$-element zero-sum mod
$r$ set with diameter at least $2s - 2 + \frac{s}{(r,s)} - 1$ in
the interval $[4s - 3, \max( last(\mathbb{Z}), last(\infty))]$.
Also, the interval $[3s - 2, 4s -4]$ must be colored entirely by
$\infty$, else condition $(i)$ would hold using elements of $R_3$.
In order to avoid a condition $(ii)$ contradiction, the interval
$[2s - 2 - \frac{s}{(r,s)} + 1, 3s - 3]$ must be colored entirely
by $R_3$.  But then the intervals $[1, s-1]$ and $[2s - 2 -
\frac{s}{(r,s)} + 1, 2s - 1 + \frac{s}{(r,s)} - 1]$ contains an
$s$-element zero-sum mod $s$set with diameter $2s -2 +
\frac{s}{(r,s)} - 1$.  This forces condition $(ii)$ and we are
finished with this case.

Finally, let $R_4$ be the only residue colored before
$first(\infty)$.  Suppose that $2s -1 + \frac{s}{(r,s)} - 1 \leq
first(R_3) \leq 3s$.  This means that there exists either an
$s$-element zero-sum mod $s$ set or an $s$-element
$\infty$-monochromatic set of diameter at most $2s -1 +
\frac{s}{(r,s)} - 1$.  Therefore, $last(R_3) \leq 5s - 1 +
\frac{s}{(r,s)} - 1$.  In addition, let $\beta = 3s - first(R_3)$.
So there are at least $s - 4 + \beta$ elements colored by $R_4$.
Also, there are at least $s - 1 + \frac{s}{(r,s)} - 1$ elements
colored by $R_4$ in the interval $[1, 2s - 1 + \frac{s}{(r,s)} -
1]$, else condition $(ii)$ holds.  If $\beta \geq 4$, then
condition $(i)$ holds using $r$ elements of residue $R_4$.  Thus,
$3s - 3 \leq first(R_3) \leq 3s$.  Suppose that we select either
an $s$-element $\infty$-monochromatic or zero-sum mod $r$ set with
minimum diameter $s - 1 + \gamma$ from $[1, 2s -1 +
\frac{s}{(r,s)} - 1]$.  Therefore, there are at most $3s -3 -
\gamma$ elements colored by $\infty$, else condition $(i)$ holds.
Furthermore, there are at most $s - 1+ \gamma$ elements colored by
$R_3$ else condition $(ii)$ holds.  This means there are at least
$2s - 2 +\frac{s}{(r,s)} - 1$ elements colored by $R_4$, and
condition $(i)$ holds.

Now suppose that $R_4$ does not have $r$ elements.  Let $3s \leq
first(R_3)$.  Then from the previous argument, we see that $R_4$
must have $r$ elements, so now suppose that $R_3$, $R_4$ and
$\infty$ each contain $r$ elements.  Therefore $last(R_4) \leq 3s
- 3$.  Now, recall that the interval $[1, 2s - 1 + \frac{s}{(r,s)}
- 1]$ contains only $R_4$ and $\infty$. Therefore, there exists an
$s$-element $\infty$-monochromatic or zero-sum mod $s$ set with
diameter at most $2s  - 2 + \frac{s}{(r,s)} - 1$.  Thus,
$last(R_3) - first(R_3) \leq 2s - 3 + \frac{s}{(r,s)} - 1$. This
means that the interval $[3s - 2, 4s - 3]$ is colored entirely by
$\infty$.  Hence the interval $[2s - 2 - \frac{s}{(r,s)} + 1, 3s -
3]$ must be colored entirely by $R_4$.  But this forces the
interval $[s  -
\frac{s}{(r,s)} + 1, 2s - 1 - \frac{s}{(r,s)} + 1]$ to be  
colored entirely by $\infty$.  However, condition $(i)$ holds by
taking $first(\infty) \cup last_1^{r-1}(\infty)$.  This completes
the case.

Let $a = first(\mathbb{Z})$ and $b = first_2(\mathbb{Z})$.  Let $y
= last_2(\mathbb{Z})$ and $z = last(\mathbb{Z})$.   Consider the
interval $S = [a, z]$.  We must satisfy the conditions of Theorem
2.3, by 1) creating an $(r-2)$-set partition with $r-3$ sets of
cardinality 2 and one set of cardinality 3, and 2) fixing two
elements that are not part of the set partition.  Let $c$ be the
residue of $\mathbb{Z}$ that colors the most integers.

\noindent \textbf{Case 3.2} Suppose there are fewer than $r$
elements colored by $c$.

There are at least $2 \lfloor \frac{5s - 1}{3} \rfloor - 3 $  $(2r
- 1 + \omega - \omega')$ elements colored by $\mathbb{Z}$, since
there can be at most $\omega' = \max(s + r - 2 - \delta, s +
\lfloor \frac{5s - 1}{3} \rfloor - 1 )$ elements colored by
$\infty$, else condition $(ii)$ holds, by taking
$first_1^s(\infty)$ and $first_{s+1}^{s+r-1}(\infty) \cup
last(\infty)$ as the $s$ and $r$-element $\infty$ monochromatic
sets, respectively.

First, suppose there are at least $2r$ elements colored by
$\mathbb{Z}$.  If this is the case, then there are at least $r +
1$ elements colored by $\mathbb{Z}$ not colored by $c$.  Fix $a$
and $z$ to be outside of the $(r - 2)$-set partition, except if in
the interval $[b,y]$ there are two colors (not including $\infty$)
each containing $r-1$ elements.  In this case, we can force $b$
and $z$ to be elements of the potential $r$-element zero-sum mod
$r$ set guaranteed in a condition of Theorem 2.3 by fixing them.
The diameter of this set is at least $3s - 3$.

Next, consider the case when there are exactly $2r - 1$ elements
colored by $\mathbb{Z}$.  then there are exactly $s + r - 1$
elements colored by $\infty$.  There must be exactly $r-1$
elements colored by a single residue, in this case $c$, else there
are at least $r+1$ elements not colored by $c$ and we can satisfy
the conditions of Theorem 2.3 in the following manner.  Fix $a$
and $z$ to be outside of the $(r-2)$-set partition except if in
the interval $[b,z]$ there are only two colors in the elements of
$\mathbb{Z}$.  If this occurs, interchange $b$ with $a$ or $z$
with $y$ in the exact manner as in the previous paragraph.  Hence,
there are exactly $r-1$ elements colored by $c$.  Let $d$ be the
residue with the second largest number of elements.  If there are
at least two elements of $\mathbb{Z}$ colored by neither $c$ or
$d$, we will construct an $(r-2)$-set partition of diameter at
least $z - b \geq 3s - 3$ or $y - a \geq 3s -3$.  Fix either $a$
and $z$, $b$ and $z$ or $b$ and $y$, and let the remaining
elements form a set partition that satisfies Theorem 2.3.  If the
desired set partition cannot be accomplished, there is only one
element of $\mathbb{Z}$ colored by neither $c$ or $d$ and we must
not fix this element or fix two elements of color $c$ or fix two
elements of color $d$, else the desired set partition does not
exist.  The minimum diameter of a set that avoids these
conditions, fixes two elements, and produces an $(r- 2)$-set
partition that we desire is  $\lfloor \frac{3r}{2} \rfloor - 2$.
By Theorem 2.3, we either obtain a $ \lfloor \frac{3r}{2} \rfloor
- 2$ diameter zero-sum mod $r$ set, or else other conditions that
will be discussed later hold.  We finish this case by describing
what happens when we do obtain the $\lfloor \frac{3r}{2} \rfloor -
2$ diameter zero-sum mod $r$ set derived from Theorem 2.3.

Suppose there exists a $\lfloor \frac{3r}{2} \rfloor - 2$ diameter
zero-sum mod $r$ set $R$.  (If $r \geq 2s -1$, condition $(i)$
follows.)  In order for there to be exactly $2r - 1$ elements
colored by $\mathbb{Z}$, $r \geq \lfloor \frac{5s -1}{3} \rfloor -
1$ since otherwise $6s - 5 + \frac{s}{(r,s)} - 1 \geq s + 3r -2$,
and condition $(ii)$ follows since there are at least $s + r$
elements colored by $\infty$.  Hence, we have constructed an
$r$-element zero-sum mod $r$ set $R$ where $diam(R) \geq \lceil
\frac{15s-1}{6} \rceil - 4$.  Hence, there are at most $\lfloor
\frac{3s + 1}{6} \rfloor$ elements colored by $\infty$, within
$[first(R), last(R)]$, else condition $(i)$ follows.  Similarly,
there are at most $2s -r - 2$ elements colored by $\mathbb{Z}$ in
the interval $[first(\infty), last(\infty)]$ else $(i)$ holds by
taking $first_1^s(\infty)$ and $first_{s+1}^{s+ r -1}(\infty) \cup
last(\infty)$ as the $s$ and $r$-element $\infty$-monochromatic
sets, respectively.

If $first(\infty) < first(R)$, either $(i)$ or $(ii)$ follows,
depending on the placement of the elements colored in $R$.  If
there exists an element of $R$ before $first_s(\infty)$, then
condition $(i)$ follows since there are more than $2s - r - 2$
elements colored by $\mathbb{Z}$ inside $[first(\infty),
last(\infty)]$, else $diam(R) \geq 3s - 2.$  Thus, $(i)$ holds by
taking $first_1^s(\infty)$ and $first_{s + 1}^{s + r - 1}(\infty)
\cup last(\infty)$ as the $s$ and $r$-element
$\infty$-monochromatic sets, respectively.  If there does not
exist an element of $R$ before $first_s(\infty)$, condition $(ii)$
holds by taking $first_1^s(\infty)$ and the $r$-element zero-sum
mod $r$ set of diameter at least $\lceil \frac{15s -1}{6} \rceil -
4$.  Therefore, $first(R) < first(\infty)$.

If $first(\infty) > 2s - 1$, then condition $(ii)$ will follow
since we can create an $s$-element zero-sum mod $s$ set from
$first_1^{2s-1}(\mathbb{Z})$, and an $r$-element
$\infty$-monochromatic set with diameter at least $2s - 2$. Hence,
$first(\infty) \leq 2s - 1$.  Within the set $R$, there are at
least $\lceil \frac{15s-1}{6} \rceil - 3$ elements colored by
$\mathbb{Z}$.  After $first(\infty)$ there can be at most $2s - r
- 2$ elements colored by $\mathbb{Z}$ in $[first(\infty),
last(\infty)]$.  Suppose $first(\infty) \leq 2s - 2$.  This
implies that there are at least $\lceil \frac{15s-1}{6} \rceil - 3
- (2s - 3) > 2s - r -2$ elements colored by $\mathbb{Z}$ in
$[first(\infty), last(\infty)]$, which produces a contradiction.
Hence $first(\infty) = 2s -2$.  If there are three or more
residues mod $s$ in $[1, 2s -2]$, then there exists an $s$-element
zero-sum mod $s$ set in $[1, 2s - 2]$ and condition $(ii)$
follows.  Hence there are two residues mod $s$ in $[1, 2s - 2]$,
each with $s - 1$ elements.  By the pigeonhole principle, one such
residue must be color $c$ mod $r$.  By fixing the first element
not colored by $last(\mathbb{Z})$ within $[1, 2s - 2]$ and
$last(\mathbb{Z})$ (or $last_2(\mathbb{Z})$ if there is only one
residue not colored $c$ or $d$), we obtain a contradiction using
Theorem 2.3 via condition $(i)$ (or other conditions that we will
discuss later hold).  Therefore, we will assume there are at least
$r$ elements colored by $c$ mod $r$.

\noindent \textbf{Case 3.3:}  Thus, there exists a color $c \in
\mathbb{Z}_r$ such that $|\Delta^{-1}(c) \cup S| \geq \lfloor
\frac{5s - 1}{3}\rfloor - 1$ $(r)$.  If $r > \lfloor \frac{5s -
1}{3}\rfloor - 1$, let $\gamma = r - (\lfloor \frac{5s -
1}{3}\rfloor - 1)$.  If $r - (\lfloor \frac{5s - 1}{3}\rfloor - 1)
$ is negative, then $\gamma = 0$.  Suppose that there are at least
$\lceil \frac{7s + 1}{3} \rceil - 2$ $(\lceil \frac{7s + 1}{3}
\rceil - 2 + \gamma)$ elements colored by $\infty$.  If there are
no elements colored by $c$ before $first_m(\infty)$ then condition
$(ii)$ holds since there exists an $s$-element
$\infty$-monochromatic set $A$ with $diam(A) \leq \lfloor \frac{5s
- 1}{3} \rfloor - 2$  $( \lfloor \frac{5s - 1}{3} \rfloor - 2
-\gamma)$ and an $r$-element $c$-monochromatic set of diameter at
least $\lfloor \frac{5s - 2}{3} \rfloor - 1$.  If $first(\infty) >
2s - 1$, we obtain a contradiction since we can create an
$s$-element zero-sum mod $s$ set from $first_1^{2s -
1}(\mathbb{Z})$ and an $r$-element $\infty$-monochromatic set with
a diameter at least $2s - 2$, namely $first_1^{r-1}(\infty) \cup
last(\infty)$.

Hence, $first(\infty) \leq 2s -1$.  Also, notice that there are at
most $\lfloor \frac{2s-1}{3} \rfloor - 1$ $(\lfloor \frac{2s-1}{3}
\rfloor - 1 - \gamma)$ holes in $[first(\infty), last(\infty)]$,
else condition $(i)$ follows in color $\infty$.  If there are
elements of $c$ both before and after $[first(\infty),
last(\infty)]$, condition $(i)$ follows since there are at least
$\lceil \frac{7s-1}{3} \rceil - 2$  $( \lceil \frac{7s-1}{3}
\rceil - 2 + \gamma)$ holes colored by $\infty$.  If there are
more than $s$ elements colored by $c$ before $first(\infty)$,
condition $(ii)$ follows, so there must be fewer than $s$ elements
colored by $c$ before $first(\infty)$.  thus, there exists at
least one element colored by $c$ after $last(\infty)$.  Hence, all
the elements colored by $c$ occur after $first(\infty)$ else
condition $(i)$ follows via $C = first_1^{r-1}(c) \cup last(c)$.
Since there are at least $s$ elements colored by $c$ after
$last(\infty)$, there can be no elements colored by $c$ before
$first_s(\infty)$, else condition $(i)$ holds since $last(c) -
first(c) \geq 3s -2$.  But then condition $(ii)$ holds by taking
an $s$-element $\infty$-monochromatic set and an $r$-element
$c$-monochromatic set, which is an $r$-element zero-sum mod $r$
set.  Thus, there are at most $\lceil \frac{7s + 1}{3} \rceil - 3$
$(\lceil \frac{7s + 1}{3} \rceil - 3 + \gamma)$ colored by
$\infty$.  This implies that there are at least $\lfloor \frac{6s-
1}{3}\rfloor - 1$  $(r + \lfloor \frac{s}{3} \rfloor)$ elements
colored by $c$.

Now suppose there are at least $\lceil \frac{6s+1}{3} \rceil - 2$
$(\lceil \frac{6s+1}{3} \rceil - 2 + \gamma)$ elements colored by
$\infty$.  We can then make the same argument as above to obtain a
contradiction.  If $first(\infty) > 2s - 1$, we obtain a
contradiction since we can create an $s$-element zero-sum mod $s$
from $first_1^{2s-1}(\mathbb{Z})$, and an $r$-element
$\infty$-monochromatic set with a diameter of at least $2s -2$.
Hence $first(\infty) \leq 2s - 1$.  Also notice that there are at
most $\lfloor \frac{3s-1}{3} \rfloor - 1$  $(\lfloor
\frac{3s-1}{3} \rfloor - 1 - \gamma)$ holes in $[first(\infty),
last(\infty)]$, else condition $(i)$ follows.  If there are
elements of $c$ both before and after $[first(\infty),
last(\infty)]$, condition $(i)$ follows since there are at least
$\lceil \frac{6s + 1 }{3} \rceil - 2 $  $(\lceil \frac{6s + 1 }{3}
\rceil - 2 + \gamma)$ holes colored by $\infty$.  If there are
more than $s$ elements colored by $c$ before $first(\infty)$,
condition $(ii)$ follows.  Hence, there must be fewer than $s$
elements colored by $c$ before $last(\infty)$.  Hence, all the
elements colored by $c$ occur after $first(\infty)$.  Since there
are at least $s$ elements colored by $c$ after $last(\infty)$,
there can be no elements colored by $c$ before $first_s(\infty)$,
else condition $(i)$ holds since $last(c)- first(c) \geq 3s -2$.
Hence, condition $(ii)$ holds by taking an $s$-element
$\infty$-monochromatic set and an $r$-element $c$-monochromatic
set, which is an $r$-element zero-sum mod $r$ set.  Therefore,
there are at most $\lceil \frac{6s + 1}{3} \rceil - 3$  $(\lceil
\frac{6s + 1}{3} \rceil - 3 + \gamma)$ elements colored by
$\infty$.  This implies there are at least $\lfloor \frac{7s -
1}{3} \rfloor - 1 $  $(r + \lfloor \frac{2s}{3} \rfloor)$ elements
colored by $c$.

Next, suppose there are at least $\lceil \frac{5s+ 2}{3} \rceil -
2$  $(\lceil \frac{5s + 2}{3} \rceil - 2 + \gamma)$ elements
colored by $\infty$.  If there are at least $\lfloor \frac{s-1}{3}
\rfloor + 1$  $(\lfloor \frac{s-1}{3} \rfloor + 1 - \gamma)$
elements colored by $\mathbb{Z}$ inside $[first(\infty),
last(\infty)]$, then we proceed in exactly the same manner as the
previous two paragraphs.  Hence, there are at most $\lfloor
\frac{s-1}{3} \rfloor$  $(\lfloor \frac{s-1}{3} \rfloor - \gamma)$
elements colored by $c$ in $[first(\infty), last(\infty)]$.  In
this case, it still follows that elements colored by $c$ can not
be both before $first(\infty)$ and after $last(\infty)$, else
condition $(i)$ holds.  Also, if the elements of $c$ are after
$last(\infty)$, condition $(ii)$ follows by selecting an
$s$-element $\infty$-monochromatic set from $[first(\infty),
last(\infty)]$ and an $r$-element $c$-monochromatic set with the
elements colored by $c$ after $last(\infty)$.  Hence, there are at
least $2s $  $(r + \lfloor \frac{s}{3} \rfloor )$ elements colored
by $c$ before $first(\infty)$.  Let $\alpha$ be the number of
elements of $\mathbb{Z}$ not colored by $c$ in the minimal
diameter $s$-element zero-sum mod $s$ set in $[first(c),
last(c)]$. Call this minimal-diameter $s$-element zero-sum mod $s$
set $K$ of diameter at most $t = s -1 + \alpha + \frac{s}{(r,s)} -
1$. If this diameter is less than or equal to $last(\infty) -
first(\infty)$, then we are done. Otherwise, we know that $t \geq
\lfloor \frac{5s-1}{3} \rfloor - 2$.  Notice that $last(\infty) -
first(\infty) \leq 2s - 3$, else condition $(ii)$ holds and there
are at most $3s - 3$ elements colored by $c$ mod $r$.  This means
there are at least $s + \alpha + \frac{s}{(r,s)} - 1$ elements
colored by $\mathbb{Z}$ but not $c$ mod $r$.  But we know that
$\alpha + \frac{s}{(r,s)} - 1 \geq \lfloor \frac{2s-1}{3} \rfloor
- 1$.  Thus, there are at least $r$ elements colored by
$\mathbb{Z}$ not colored by $c$, and we have already considered
this case.  Thus, there are at least $\lfloor \frac{8s-1}{3}
\rfloor - 1$  $( \lfloor \frac{8s-1}{3} \rfloor - 1 + \gamma)$
elements colored by $c$ and there are at most $\lceil \frac{5s +
1}{3} \rceil - 3$  $( \lceil \frac{5s + 1}{3} \rceil - 3+ \gamma <
r)$ elements colored by $\infty$.  This contradicts our assumption
that $\Delta^{-1}(\infty)| \geq r$. Thus, we are finished with
this case.

Therefore, we can now apply Theorem 2.3 to observe that either the
lemma is satisfied or all but at most $a - 2$ of the elements
colored by $\mathbb{Z}$ are colored by elements from the same
coset $(a\mathbb{Z}_m + \alpha)$ of $\mathbb{Z}_m$.  In the latter
case, Theorem 2.2 implies that any subset of cardinality $(m +
\frac{m}{a} + a - 3 ) \leq \lceil \frac{3s}{2} \rceil - 1$ must
contain an $m$-element zero-sum mod $m$ set.

If $r \geq 2s -2$, condition $(i)$ follows with the $r$-element
zero-sum mod $r$ set $G'$ constructed from the set $G =
first_1^{\lceil \frac{r}{2} \rceil}(\mathbb{Z}) \cup
last_1^{r-1}(\mathbb{Z})$.  If there are at least $\lceil
\frac{r}{2} \rceil$ elements colored by $\mathbb{Z}$ before
$first(\infty)$, then condition $(i)$ or $(ii)$ follows. Condition
$(i)$ follows if $last_{\lceil \frac{r}{2} \rceil }(\mathbb{Z})
\geq 3s - 3$, since $G'$ has a diameter of at least $3s - 3$.
However, if this is not the case, then there are at most $\lceil
\frac{r}{2} \rceil - 1$ integers colored by $\mathbb{Z}$ in the
interval $[3s - 2, \max(last(\mathbb{Z}), last(\infty))]$.  Hence,
$last(\infty) \geq 6s - 5 + \frac{s}{(r,s)} - 1 - \lceil
\frac{r}{2} \rceil$  $(last(\infty) \geq 2r - 1 + \omega  - \lceil
\frac{r}{2} \rceil)$, and therefore $first(\infty) \geq 3s - 2 -
\lceil \frac{r}{2} \rceil$ $(3r - 2s - 2 - \delta - \lceil
\frac{r}{2} \rceil)$. Hence, the first $\lceil \frac{3s}{2}
\rceil$ elements of $[1, \max(last(\mathbb{Z}), last(\infty))]$
are colored entirely by $\mathbb{Z}$, and thus there exists an
$s$-element zero-sum mod $s$ set $H$ with $diam(H) \leq \lceil
\frac{3s}{2} \rceil - 2$. Hence, condition $(ii)$ holds.
Therefore, $first(\infty) \leq \lceil \frac{r}{2} \rceil$.  Thus,
the interval $[3s + \lceil \frac{r}{2} \rceil - 2 ,
\max(last(\mathbb{Z}), last(\infty))]$ is colored entirely by
$\mathbb{Z}$.  Suppose there are fewer than $s$ elements colored
by $\mathbb{Z}$ between $[first(\infty), last(\infty)]$.  Let $A_1
= [1, first(\infty) - 1], A_2 = last_1^{\lceil \frac{r}{2} \rceil
- first(\infty) + 1}(\mathbb{Z}, [last(\infty), 4s - 3])$
$(last_1^{\lceil \frac{r}{2} \rceil - first(\infty) +
1}(\mathbb{Z}, [last(\infty)$ $, \lceil \frac{2s}{3} \rceil + 2r -
3 - \delta]))$, and $A_3 = [6s - r - 5, \max(last(\mathbb{Z}),
last(\infty))]$  $( [r - 1 + \omega, \max(last(\mathbb{Z}),
last(\infty))])$.  In this case, we are guaranteed that $A = A_1
\cup A_2 \cup A_3$ contains an $r$-element zero-sum mod $r$ set.
Condition $(i)$ holds if the set consists of elements from $A_1$
and $A_3$, or $A_1, A_2$ and $A_3$.  An $r$-element zero-sum mod
$r$ set can not be constructed using only $A_1$ and $A_2$.  If an
$r$-element zero-sum mod $r$ set appears in $A_2$ and $A_3$ only,
then condition $(ii)$ follows since there exists an $s$-element
$\infty$-monochromatic set $B$ with $diam(B) \leq 2s -2$ using the
first $s$ elements colored by $\infty$, and there exists an
$r$-element zero-sum mod $r$ set $B'$ with $diam(B') \geq 2s - 1$.
Thus there must be at least $s$ elements colored by $\mathbb{Z}$
in $[first(\infty), last(\infty)]$.  This implies there are at
most $2s - 3$ elements colored by $\infty$ else condition $(i)$
holds.  Thus, there are at least $4s - 2 (2r - 2s + 2 + \omega)$
elements colored by $\mathbb{Z}$.  In this case, the $r$-element
zero-sum mod $r$ set $A'$ as described before, has a diameter of
at least $3s - 3$, and thus condition $(i)$ follows. We have a
contradiction, and the proof is complete.  \qed
\end{proof}

\begin{theorem}
Let $r$ and $s$ be positive integers where $r \geq s \geq 3$. If
$r \leq \lfloor \frac{5s-1}{3} \rfloor - 1$, then $f(s, r,
\{\infty\} \cup \mathbb{Z}) = 9s - 7 + \frac{s}{(r,s)} - 1$.
\end{theorem}
\begin{proof}
The coloring $\Delta : [1, 9s - 8 + \frac{s}{(r,s)} - 1]
\longrightarrow \{\infty\} \cup \mathbb{Z} $ given by the string:
\begin{align*}
10^{s-1}1^{s-2}\infty^{s-1}10^{s-1}2^{\frac{s}{(r,s)} - 1}0^{s-1}1^{2s - 2}\infty^{2s - 2} & \textnormal{   mod } s \\
10^{s-1}1^{s-2}\infty^{s-1}10^{2s - 2 + \frac{s}{(r,s)} - 1}1^{2s
- 2}\infty^{2s - 2} & \textnormal{ mod } r
\end{align*}
implies that $f(s,r,\{\infty\} \cup \mathbb{Z}) > 9s - 8 +
\frac{s}{(r,s)} - 1$. Next we show that $f(s,r, \{\infty\} \cup
\mathbb{Z}) \leq 9s - 7 + \frac{s}{(r,s)} - 1.$  Let $\Delta : [1,
9s - 7 + \frac{s}{(r,s)} - 1] \longrightarrow \{\infty\} \cup
\mathbb{Z})$ be an arbitrary coloring. By Theorem 1.1, the
interval $[1, 3s - 2]$ contains an $s$-element zero-sum mod $s$
subset $S$ such that $diam(S) \leq 3s - 3$. Shifting the interval
$[3s - 1, 9s - 7 + \frac{s}{(r,s)} - 1]$ to the interval $[1, 6s -
5 + \frac{s}{(r,s)} - 1]$ and applying Lemma 4.7 completes the
proof for $r \geq s \geq 3$. \qed
\end{proof}

\begin{theorem}
Let $r$ and $s$ be positive integers where $r \geq s \geq 3.$  If
$ \lfloor \frac{5s-1}{3} \rfloor - 1 \leq r \leq 3s -3$, then
$f(s, r, \{\infty\} \cup \mathbb{Z}) = 3s + 2r - 3 + \min(3s - 3,
s + r - 1 + \frac{s}{(r,s)} - 1)$.
\end{theorem}
\begin{proof}
The coloring $\Delta : [1, 3s + 2r - 4 + \min(3s - 3, s + r - 1 +
\frac{s}{(r,s)} - 1)] \longrightarrow \{\infty\} \cup \mathbb{Z} $
given by the string:
\begin{align*}
10^{s-1}1^{s-2}\infty^{s-1}10^{s-1}2^{\frac{s}{(r,s)} - 1}0^{\min(r, 2s - 1 - \frac{s}{(r,s)})}1^{r - 1}\infty^{r - 1} & \textnormal{   mod } s \\
10^{s-1}1^{s-2}\infty^{s-1}10^{\min(3s - 3, s + r - 1 +
\frac{s}{(r,s)} - 1)}1^{r-1}\infty^{r-1} & \textnormal{ mod } r
\end{align*}
implies that $f(s,r,\{\infty\} \cup \mathbb{Z}) > 3s + 2r - 4 +
\min(3s - 3, s + r - 1 + \frac{s}{(r,s)} - 1)$. Next we show that
$f(s,r, \{\infty\} \cup \mathbb{Z}) \leq 3s + 2r - 3 + \min(3s -
3, s + r - 1 + \frac{s}{(r,s)} - 1).$  Let $\Delta : [1, 3s + 2r -
3 + \min(3s - 3, s + r - 1 + \frac{s}{(r,s)} - 1)] \longrightarrow
\{\infty\} \cup \mathbb{Z})$ be an arbitrary coloring. By Theorem
1.1, the interval $[1, 3s - 2]$ contains an $s$-element zero-sum
mod $s$ subset $S$ such that $diam(S) \leq 3s - 3$. Shifting the
interval $[3s - 1, 3s + 2r - 3 + \min(3s - 3, s + r - 1 +
\frac{s}{(r,s)} - 1)]$ to the interval $[1,  2r - 1 + \min(3s - 3,
s + r - 1 + \frac{s}{(r,s)} - 1)]$ and applying Lemma 4.7
completes the proof for $r \geq s \geq 3$. \qed
\end{proof}
Notice that for these two cases, we have a generalization in the
sense of Erd\H{o}s-Ginzburg-Ziv if and only if $(r,s) = 1$ or $3s
- 3 \geq s + r - 1 + \frac{s}{(r,s)} - 1$.  As we can see in the
next theorem, for $r
> 3s - 3$ we always obtain a generalization in the sense of
Erd\H{o}s-Ginzburg-Ziv.

\begin{theorem}
Let $r$ and $s$ be positive integers where $r \geq s \geq 2$.  If
$r > 3s - 3$, then $f(s, r, \{\infty\} \cup \mathbb{Z}) = 3s + 3r
- 4$.
\end{theorem}
\begin{proof}
The coloring $\Delta : [1, 3s + 3r - 5] \longrightarrow \{\infty\}
\cup \mathbb{Z} $ given by the following string, which is the same
mod $s$ or mod $r$:
\begin{align*}
\infty1^{s-1}\infty^{s-2}0^{s-1}\infty1^{r-1}\infty^{r-1}0^{r-1}
\end{align*}
implies that $f(s,r, \{\infty\} \cup \mathbb{Z}) \geq 3s + 3r -
4$.  By the pigeonhole principle, it follows that the interval
$[1, 3s - 2]$ contains an $s$-element zero-sum mod $s$ or
$\infty$-monochromatic subset $S_1$ with $diam(S_1) \leq 3s - 3$,
and that the interval $[3s - 1, 3s + 3r - 4]$ contains an
$r$-element zero-sum mod $r$ or $\infty$-monochromatic subset
$S_2$ with $diam(S_2) \geq 3s - 3$.  Hence $S_1$ and $S_2$
complete the proof.  \qed
\end{proof}

\begin{acknowledgement}
The author would like to thank Arie Bialostocki for supervising
this research and to David Grynkiewicz for his helpful
communications.
\end{acknowledgement}

\end{document}